\documentclass[a4paper,12pt]{article}
\usepackage[english]{babel} 
\usepackage{graphicx}
\usepackage{amsmath}
\usepackage{amsfonts}
\usepackage[square,comma,numbers,sort&compress]{natbib}
\usepackage{ifthen}
\usepackage{curves}
\usepackage{times,fancyhdr}

\newtheorem{thrm}{Theorem}[section]
\newtheorem{dfntn}[thrm]{Definition}
\newtheorem{lmm}[thrm]{Lemma}
\newtheorem{crllr}[thrm]{Corollary}

\def\C{{\mathbb C}}              % Complex numbers
\def\R{{\mathbb R}}              % Real numbers
\def\Z{{\mathbb Z}}              % Integer numbers

\def\Pr{{\mathbb P}}             % Probability of an event
\def\E{{\mathbb E}}              % Expectation 
\def\T{{\mathbb T}}              % Stopping time
\def\I{{\mathbb I}}              % Indicator function

\def\M{{\mathcal M}}             % Dirichlet-to-Neumann operator
\def\B{{\mathcal B}}             % Borel sigma-algebra
\def\N{{\mathcal N}}             % Number of reflections

\def\l{\ell}                     % Minimal cut-off or local time

\def\ve{\varepsilon }            % Reflection probability

\def\pa{{\partial\Omega}}        % Boundary of the domain

\def\V{{\bf V}}                  % Eigenfunction of the operator M
\def\P{{\bf P}}                  % Distribution of hitting or absorption probabilities

\def\TT{{\mathcal T}}

\def\sign{\textrm{sign}}

\def\<{< \hspace*{-0.5mm}}
\def\>{\hspace*{-0.5mm} >}

\def\Prob{{\mathcal P}}          % Probability of the absorption on a disk
\def\K{{\mathcal K}}             % Gaussian error function

\def\LL{{L^2}}                       % L^2(\pa) space 

\begin{document}
\title{  \vspace*{-5mm}
Partially Reflected Brownian Motion: \\A Stochastic Approach to Transport Phenomena
\footnote{ This article partially reproduces the chapter which 
has been written by the author for the volume ``Focus on 
Probability Theory'', and it should be referenced as 
D.~S.~Grebenkov, in {\it Focus on Probability Theory}, 
Ed. L.~R.~Velle, pp. 135-169 (Nova Science Publishers, 2006). 
The bibligraphic reference (ISBN) is 1-59454-474-3. 
Further information about this volume can be found
on https://www.novapublishers.com/catalog/ } }
\author{\bfseries\itshape Denis S.~Grebenkov\thanks{E-mail address: denis.grebenkov@polytechnique.edu}\\
\small Laboratoire de Physique de la Mati\`ere Condens\'ee,\\ 
\small C.N.R.S. -- Ecole Polytechnique 91128 Palaiseau Cedex, France}
\date{\small \it Received: August 2004; \hskip 10mm  Published: September 2006}
\maketitle
%\thispagestyle{empty}
%\setcounter{page}{1}
% ------------ [First Page Running Head] -------------------------------------------------
%\thispagestyle{fancy}
%\fancyhead{}
%\fancyhead[L]{In: Book Title \\ 
%Editor: Editor Name, pp. {\thepage-\pageref{lastpage-01}}}
%\fancyhead[R]{ISBN 0000000000  \\
%\copyright~2006 Nova Science Publishers, Inc.}
%\fancyfoot{}
%\renewcommand{\headrulewidth}{0pt}
%-------------------------------------------------------------------------------

\begin{abstract}
Transport phenomena are ubiquitous in nature and known to be important
for various scientific domains. Examples can be found in physics,
electrochemistry, heterogeneous catalysis, physiology, etc.  To obtain
new information about diffusive or Laplacian transport towards a
semi-permeable or resistive interface, one can study the random
trajectories of diffusing particles modeled, in a first approximation,
by the partially reflected Brownian motion.  This stochastic process
turns out to be a convenient mathematical foundation for discrete,
semi-continuous and continuous theoretical descriptions of diffusive
transport.

This paper presents an overview of these topics with a special
emphasis on the close relation between stochastic processes with
partial reflections and Laplacian transport phenomena. We give
selected examples of these phenomena followed by a brief introduction
to the partially reflected Brownian motion and related probabilistic
topics (e.g., local time process and spread harmonic measure). A
particular attention is paid to the use of the Dirichlet-to-Neumann
operator. Some practical consequences and further perspectives are
discussed.
\end{abstract}

%\subjclass{ 60J60, 31C20, 60G50, 80A20 }

% 60J60  Diffusional processes
% 31C20  Discrete potential theory and numerical methods
% 60G50  Sum of independent random variables; random walks
% 80A20  Heat and mass transfer, heat flow


\section*{ Introduction }

An erratic motion of pollens of Clarkia (primrose family) discovered
by Robert Brown in 1827 and quantitatively described by Albert
Einstein in 1905 gave a substantial impact for developing mathematical
theory of stochastic processes, an important branch of modern
mathematics. Supported by rigorous mathematical foundations, the
Brownian motion and related stochastic processes found numerous
applications in different scientific domains, from theoretical physics
to biology and economics. To study the transport of species diffusing
from a remote source towards and across semi-permeable or resistive
interface (e.g., cellular membrane), one can employ either an averaged
description in terms of an appropriate boundary value problem for the
concentration of species, or stochastic analysis of their random
trajectories. In the first case, a finite permeability (reactivity,
resistivity, etc.) of the interface leads to the mixed or Fourier
boundary condition, while in the second case it can be modeled as
partial reflections on the boundary.  Physical or chemical processes
governed by the Laplace equation (stationary diffusion) with mixed
boundary condition are generally called {\it Laplacian transport
phenomena}. Their examples are found in physiology (oxygen diffusion
towards and across alveolar membranes), in electrochemistry (electric
transport in electrolytic cells), in heterogeneous catalysis
(diffusion of reactive molecules towards catalytic surfaces), in
nuclear magnetic resonance (diffusion of spins in confining porous
media). Studying random trajectories of diffusing species, one can
extract a subtle information about the system in question.  However,
the theoretical or numerical analysis of these phenomena is in general
complicated by an irregular geometry of the interface (e.g.,
microroughness of metallic electrodes, see Section~\ref{sec:PTL}).

In this paper, we focus on a particular stochastic process, called
{\it partially reflected Brownian motion} (PRBM), and its application
to study Laplacian transport phenomena. Our main purpose is to capture
the attention to this interesting process itself, and its use for
understanding the influence of a geometrical irregularity of the
interface on transport properties. Bearing in mind the particular role
of the geometry, we would like to ``bridge'' theoretical, numerical
and experimental studies of Laplacian transport phenomena on the one
side, and powerful mathematical methods of stochastic analysis on the
other side. Since the extensive literature existing on both topics is
generally difficult to get through for non-specialists, we prefer to
use a descriptive style of writing in order to give the whole vision
of the problem, without specifying particular details which can be
found anywhere else (e.g., see references at the end of this paper).

% ------------ [Running Heads - pagina 2] -------------------------------------------------
\pagestyle{fancy}
\fancyhead{}
\fancyhead[EC]{Denis S.~Grebenkov}
\fancyhead[EL,OR]{\thepage}
\fancyhead[OC]{Partially Reflected Brownian Motion...}
\fancyfoot{}
\renewcommand\headrulewidth{0.5pt} 
%-------------------------------------------------------------------------------

In the first section, we present three examples of Laplacian transport
phenomena in different scientific fields. Their mathematical
description by the mixed boundary value problem opens encouraging
possibilities to apply powerful tools of potential theory, variational
analysis and probability theory. The second section is devoted to
remind some basic definitions of stochastic process theory: stopping
times, reflected Brownian motion, local time process, harmonic
measure, etc. In the third section, we introduce the partially
reflected Brownian motion and show its properties for a planar
surface. An important relation to the Dirichlet-to-Neumann operator is
revealed and then illustrated by several examples. The last section
presents different stochastic descriptions of Laplacian transport
phenomena: a recently developed continuous approach and two other
methods.  In the conclusion, we summarize the essential issues of the
paper.

\section{ Laplacian Transport Phenomena }
\label{sec:PTL}

The transport of species between two distinct ``regions'' separated by
an interface occurs in various biological systems: water and minerals
are pumped by roots from the earth, ions and biological species
penetrate through cellular membranes, oxygen molecules diffuse towards
and pass through alveolar ducts, and so on. Transport processes are
relevant for many other scientific domains, for example, heterogeneous
catalysis and electrochemistry.  In this section, we shall give three%
\footnote{
Diffusive NMR phenomena present another important example when the
transport properties are considerably affected by irregular geometry.
In this paper, we do not discuss this case since it has been recently
reviewed in a separate paper \cite{Grebenkov06}. }
%footnote
important examples of the particular transport process, called {\it
Laplacian} or {\it diffusive transport}.

\subsection{ Stationary Diffusion across Semi-permeable Membranes }

\begin{figure}
\begin{center}
\includegraphics[width=65mm]{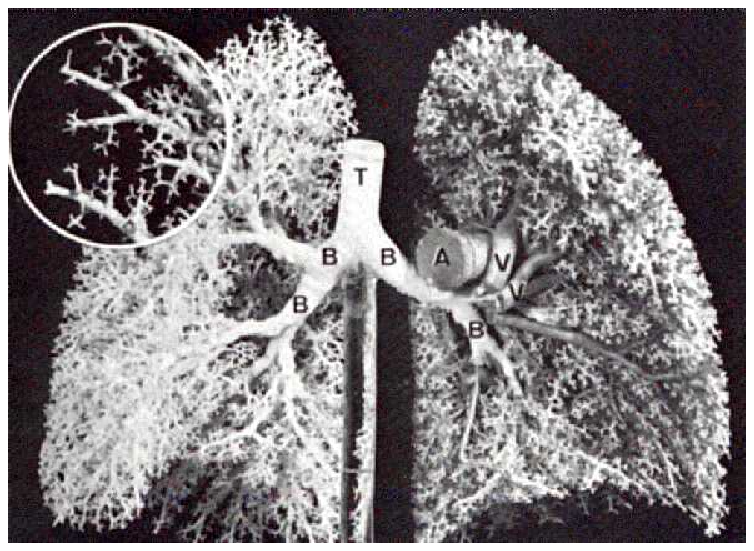}  
\hskip 6mm  
\includegraphics[width=64mm]{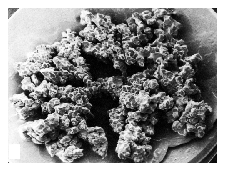}
\end{center}
\caption{ 
On the left, a cast of human lungs; on the right, a cast of pulmonary
acinus \cite{Weibel,Haefeli88} (by E.~Weibel).}
\label{fig:lung}
\end{figure}

Let us begin by considering the respiration process of
mammals. Inbreathing a fresh air, one makes it flow from the mouth to
the dichotomic bronchial tree of the lungs (Fig.~\ref{fig:lung}). For
humans, first fifteen generations of this tree serve for convectional
transport of the air towards pulmonary acini, terminal gas exchange
units \cite{Weibel}. A gradual increase of the total cross-section of
bronchiae leads to a decrease of air velocity. At the entrance of the
acinus, it becomes lower than the characteristic diffusion velocity
\cite{Mauroy04}. As a consequence, one can describe the gas exchange
inside the acinus as stationary diffusion of oxygen molecules in air
from the entrance (``source'' with constant concentration $C_0$ during
one cycle of respiration) to the alveolar membranes \cite{Sapoval01b}.
In the bulk, the flux density is proportional to the gradient of
concentration (Fick's law), ${\bf J}=-D\nabla C$, where $D$ is the
diffusion coefficient. The mass conservation law, written locally as
$\textrm{div}~ {\bf J}=0$, leads to the Laplace equation $\Delta C=0$
in the bulk. The flux density {\it towards} the interface is simply
$J_n=D~\partial C/\partial n$, where the normal derivative $\partial
/\partial n$ is directed to the bulk. Arrived to the alveolar
membrane, oxygen molecules can penetrate across the boundary for
further absorption in blood, or to be ``bounced'' on it and to
continue the motion.  The ``proportion'' of absorbed and reflected
molecules can be characterized by {\it permeability} $W$ varying from
$0$ (perfectly reflecting boundary) to infinity (perfectly absorbing
boundary). In this description, the flux density {\it across} the
alveolar membrane is proportional to the concentration, $J_n =WC$.
Equating these two densities on the alveolar membrane, one gets a
mixed boundary condition, $D(\partial C/\partial n)=WC$, called also
{\it Fourier} or {\it Robin} boundary condition. Resuming these
relations, one provides the following mathematical description for the
diffusion regime of human or, in general, mammalian respiration:
\begin{eqnarray}                                                                                   
\label{eq:problem_C1}
\Delta C = 0 ~~ &&  \textrm{in the bulk} \\                                                           
\label{eq:problem_C2}
C = C_0      &&  \textrm{on the source} \\                                                         
\label{eq:problem_C3}
\left[I-\Lambda \frac{\partial }{\partial n}\right] C = 0 ~~
&& \textrm{on the alveolar membrane} 
\end{eqnarray}
where the underlying physics and physiology are characterized by a
single parameter $\Lambda =D/W$, which is homogeneous to a length
($I$ stands for the identity operator).  Note also that the dependence
on constant $C_0$ is irrelevant.  In what follows, we address to this
``classical'' boundary value problem. The essential complication
resides in a very irregular geometry of the pulmonary acinus, which
presents a branched structure of eight generations (for humans),
``sticked'' by alveolar ducts (Fig.~\ref{fig:lung}). For small
$\Lambda$, only a minor part of the boundary is involved to the
transport process (so-called {\it Dirichlet active zone}), whereas the
flux across the rest of the boundary is almost zero (this effect is
called {\it diffusional screening}
\cite{Sapoval01b,Sapoval02,Sapoval02a,Felici03,Felici04,Grebenkov05a}).
With an increase of $\Lambda$, larger and larger part of the boundary
becomes active.  As a result, the efficiency of human lungs depends on
the parameter $\Lambda$ in a nontrivial manner that implies different
physiological consequences \cite{Felici05}. The trajectory of a chosen
oxygen molecule can be seen as Brownian motion from the source towards
the alveolar membrane, with multiple bounces on the boundary and final
absorption.  This is in fact what we call the partially reflected
Brownian motion (Section~\ref{sec:PRBM}). A profound study of the
interplay between the irregular geometry of the acinus and the erratic
random motion of oxygen molecules inside it should help to better
understand physiological functioning of human lungs.

\subsection{ Heterogeneous Catalysis }

\begin{figure}
\begin{center}
\includegraphics[width=45mm]{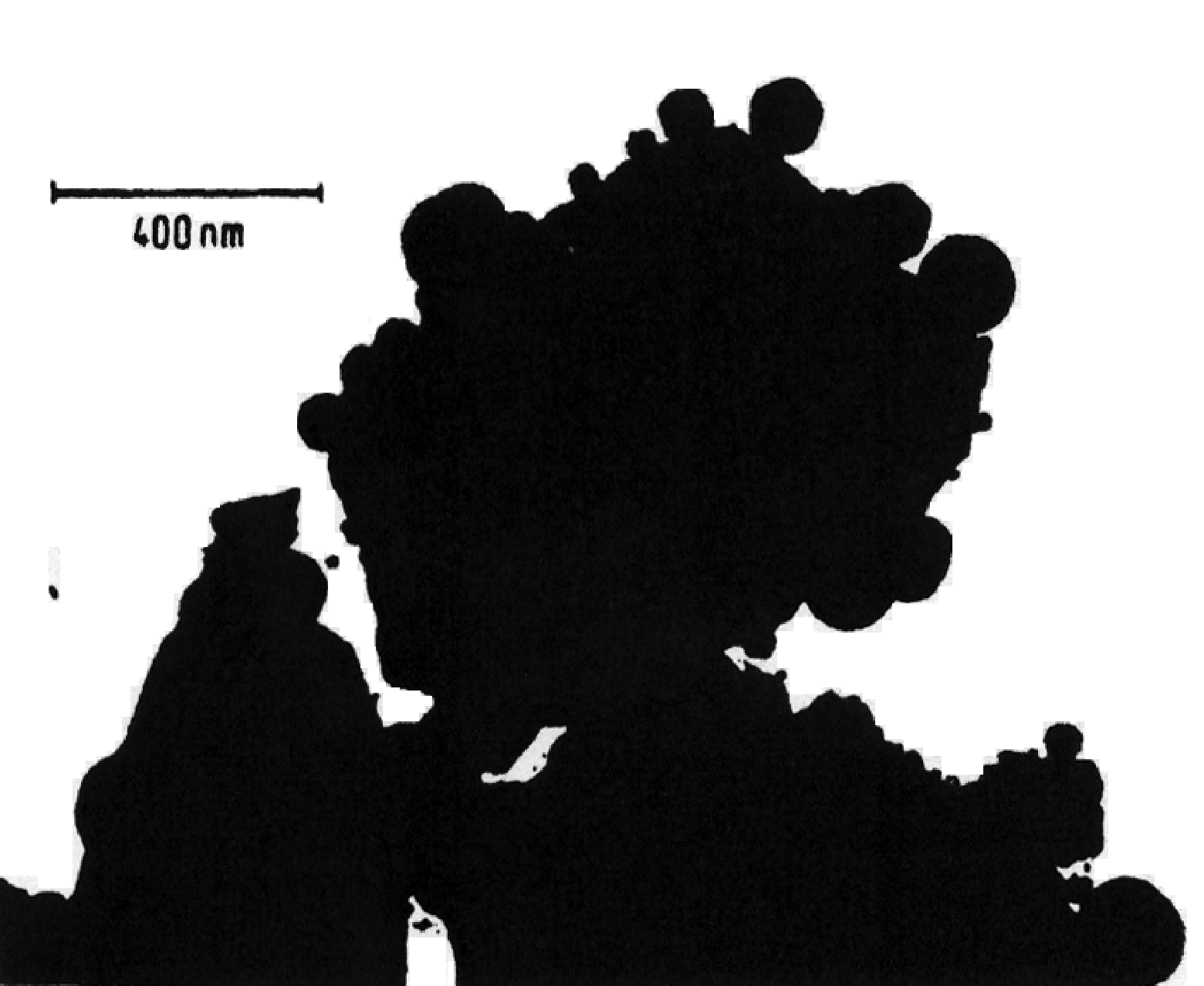}  \hskip 4mm
\includegraphics[width=45mm]{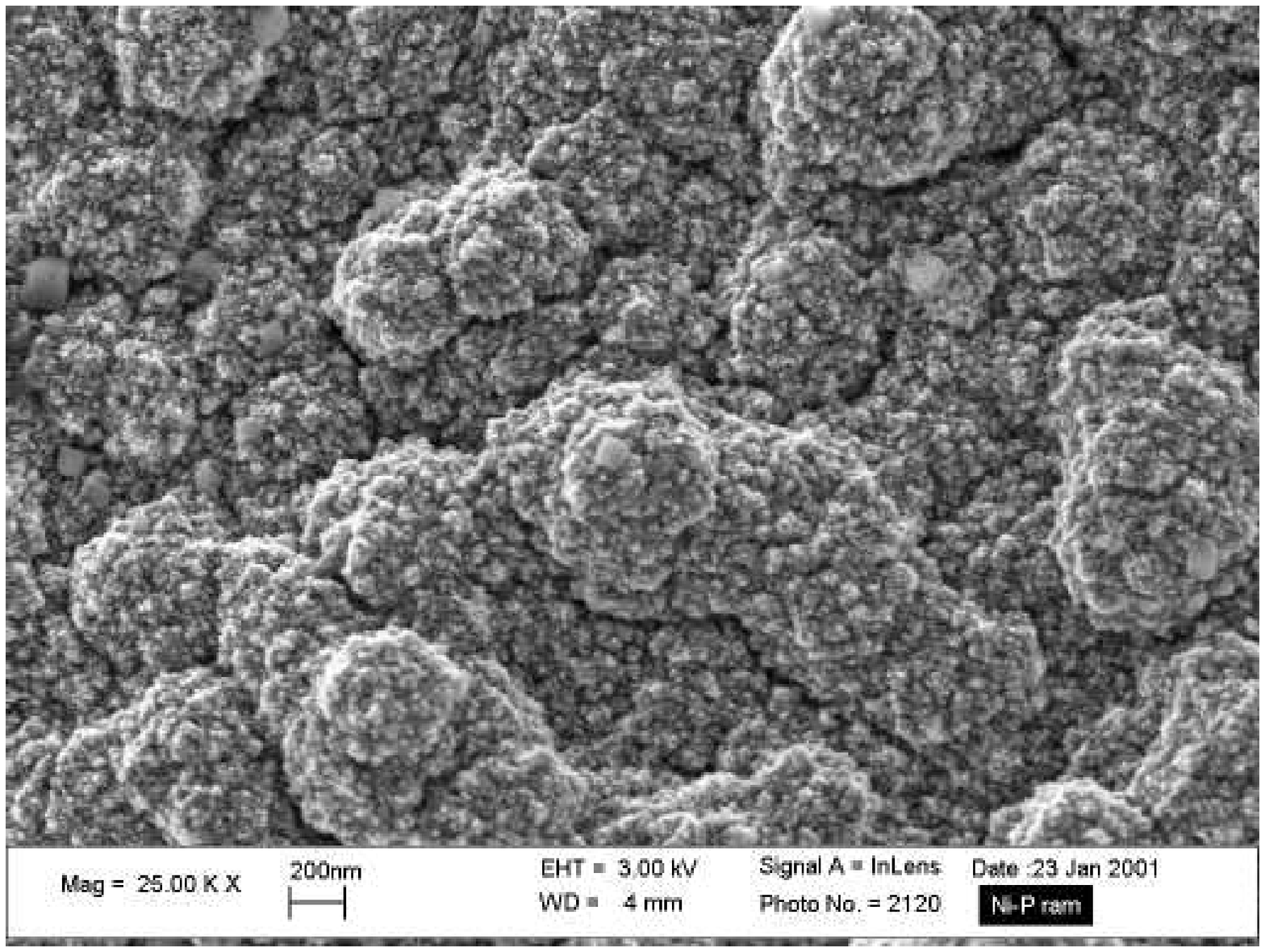}  \hskip 4mm
\includegraphics[width=35mm]{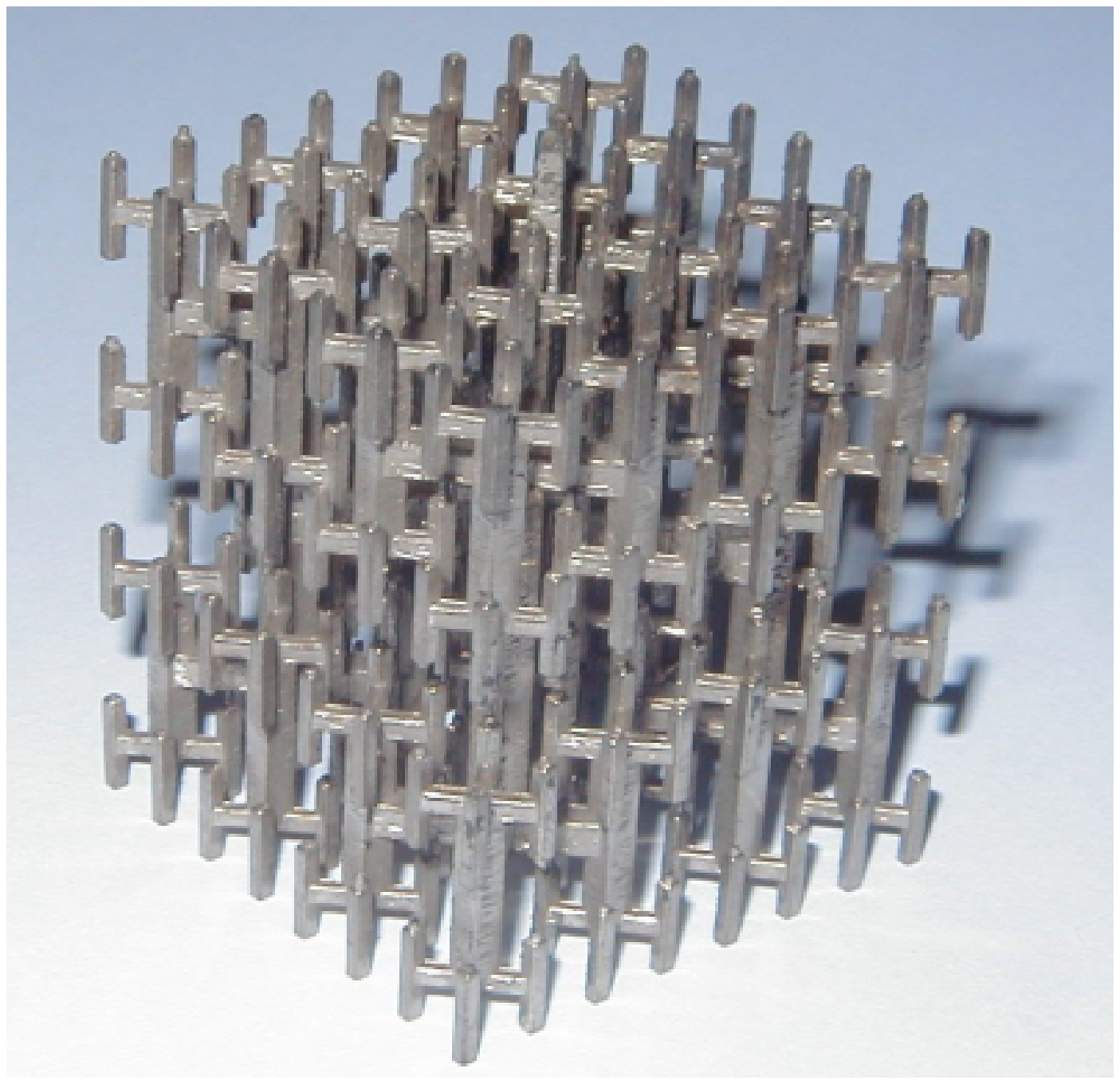}
\end{center}
\caption{ 
On the left, an example of an irregular catalytic surface (by
J.~S.~Andrade jr.); at the center, photo of a rough metallic surface
of nickel electrode (by E.~Chassaing); on the right, photo of an
irregular metallic electrode used to study the Laplacian transport
phenomena experimentally (by B.~Sapoval).}
\label{fig:catalyst}
\end{figure}

A similar description can be brought to the molecular regime of
heterogeneous catalysis omnipresent in petrochemistry. One considers
reactive molecules $A$ injected into a solvent and then diffusing
towards a catalyst. Hitting the catalytic surface, they can be
transformed into other molecules $A^*$ (with a finite reaction rate
$K$), or to be bounced for further diffusion in the bulk. The new
molecules $A^*$, collected by appropriate physical or chemical
technique, do not further contribute to the transport process.
Assuming the presence of a remote source of reactive molecules $A$,
one can model, in a first approximation%
\footnote{ 
This description is probably too simplified in order to model the
heterogeneous catalysis {\it quantitatively}. First, the presence of
molecules $A^*$ near the catalyst may ``obstruct'' the access to the
catalytic surface.  Second, parasite reactions happen on the boundary
that implies a progressive deactivation of the catalyst. Consequently,
the reactivity $K$ becomes dependent on the spatial position on the
catalytic surface, leading to an inhomogeneous boundary condition.
Finally, the molecular diffusion can be applied only if the mean free
path of reactive molecules is much lower than the geometrical features
of the catalyst (in the opposite case, one deals with Knudsen
diffusion \cite{Sahimi90,Coppens99,Malek01,Andrade03}). Nevertheless,
the simple description (\ref{eq:problem_C1}--\ref{eq:problem_C3})
permits to take into account many important features related to the
catalytic process.},
%footnote
the heterogeneous catalysis by the mixed boundary value problem
(\ref{eq:problem_C1}--\ref{eq:problem_C3}) with a characteristic
length $\Lambda =D/K$ \cite{Sapoval01,Andrade01,Andrade03,Filoche05}.
The keynote of this similitude is related to the fact that each
reactive molecule arrived onto the boundary terminates its motion
after a number of successive reflections. The mechanism leading to its
termination is different: for the oxygen diffusion, the molecules are
absorbed by the alveolar membrane and transferred to the blood, while
for the heterogeneous catalysis, the reactive molecules are
transformed by chemical reaction into other molecules which do not
further participate to the process. Since the overall production of
new molecules $A^*$ depends on the total surface area of the catalytic
surface, one tries to design catalysts with the largest possible
surface (for given volume), realizing porous and very irregular
boundaries (Fig.~\ref{fig:catalyst}).  As a consequence, the
diffusional screening becomes important to understand numerous
industrial processes in petrochemistry. Since random trajectories of
reactive molecules correspond to the partially reflected Brownian
motion, its study may allow a design of more efficient catalysts.

\subsection{ Electric Transport in Electrochemistry } 

The other example of Laplacian transport phenomena can be found in
electrochemistry: the electric current between two metallic electrodes
into an electrolyte is described by the same boundary value problem.
Indeed, the electric potential $V$ obeys the Laplace equation in the
bulk since the electrolyte is locally neutral.  Taking one electrode
of very low resistance ({\it counter-electrode}), one writes the
corresponding boundary condition as $V=V_0$, where $V_0$ is the
applied tension.  For the other electrode of surface resistance $r$
({\it working electrode}), one obtains the mixed boundary condition by
equating the volume current density $-\rho^{-1}\nabla V$ ($\rho$ is
the electrolyte resistivity) and the surface current density $V/r$:
$\Lambda ~\partial V/\partial n=V$, where $\Lambda =r/\rho$ is again
the physical length of the problem. The similar description can be
brought even in the case of an alternative tension
\cite{Sapoval94,Sapoval96b}. 

For electric transport, one cannot associate directly the mixed
boundary value problem with the partially reflected Brownian motion
since there is no diffusing particle. From this point of view, the
electrochemical problem has only a formal analogy with two previous
examples. At the same time, the electrochemistry is an appropriate
domain to study experimentally the influence of the irregular geometry
on the (average) transport properties. Taking metallic electrodes of
different shapes with micro- or macro-roughness (e.g., see
Fig.~\ref{fig:catalyst}), one can directly measure the {\it
spectroscopic impedance} or {\it admittance} (see below). These
characteristics are equivalent to the total flux across the boundary
for diffusional problems \cite{Filoche99,Grebenkov03}.  The
observation of anomalous impedance behavior \cite{Wolff26} had
provoked numerous theoretical, numerical and experimental studies of
the role of a geometrical irregularity in Laplacian transport
phenomena
\cite{deLevie65,Nyikos85,Halsey87,Halsey91a,Halsey91b,Halsey92,%
Sapoval87,Sapoval88,Sapoval89,Filoche97,Sapoval99,Filoche00,%
Grebenkov05,Grebenkov05b,Levitz06}
(for more information, see \cite{Grebenkov} and references therein).

\subsection{ Discrete and Semi-continuous Approaches }

Among different theoretical approaches developed to study Laplacian
transport phenomena, we have to mention the double layer theory of
Halsey and Leibig \cite{Halsey87,Halsey91a,Halsey91b,Halsey92} and the
formalism of the Brownian self-transport operator proposed by Filoche
and Sapoval \cite{Filoche99}. In Section~\ref{sec:approaches_Halsey},
we shall show how the original Green function description by Halsey
and Leibig can be related to the Brownian motion reflected with jump
(this stochastic reformulation will be referred to as
``semi-continuous'' approach). In turn, Filoche and Sapoval considered
lattice random walks with partial reflections to derive a spectral
representation for the macroscopic response of an irregular interface
(see Section~\ref{sec:approaches_Filoche}, where this formalism is
referred to as ``discrete'' approach).  Although both methods
accurately describe Laplacian transport phenomena (e.g., they give an
explicit formula for the total flux across the boundary), their major
inconvenience resides in the dependence on an artificial length scale:
jump distance $a$ for the semi-continuous approach and lattice
parameter $a$ for the discrete approach. A physical intuition suggests
that, if these descriptions are correct, there should exist a well
defined continuous limit as $a$ tends to $0$. Certain substantial
arguments to justify the existence of this limit were brought in
\cite{Grebenkov} (and they will be strengthened in this paper), but a
rigorous mathematical proof is still required. To overcome this
difficulty, a new theoretical approach has been recently developed in
\cite{Grebenkov06b}. We shall call it ``continuous'' since it is
tightly related to a continuous stochastic process, namely, the
partially reflected Brownian motion. This approach will integrate the
advantages of the previous ones, being a mathematical foundation for
understanding Laplacian transport phenomena. We shall return to these
questions in Section~\ref{sec:approaches}.

\section{ Basic Definitions }

In this section, we recall the basic definitions related to the
Brownian motion and reflected Brownian motion that can be found in
extensive literature, e.g.,
\cite{Freidlin,Borodin,Ito,Levy,Revuz,Bass,Port}.  The familiar reader
may pass over this section.

\subsection{ Brownian Motion and Dirichlet Boundary Value Problem }

The Brownian motion can be defined in different ways \cite{Freidlin}.
Throughout this paper, we use the following definition.

\begin{dfntn}
A stochastic process $W_t$ ($t\geq 0$) defined on the chosen
probabilistic space is called {\it one-dimensional Brownian motion}
(or Wiener process) started from the origin, if
\begin{itemize}
\item
its trajectories are continuous almost surely (with probability $1$);

\item
it starts from the origin almost surely, $\Pr\{ W_0=0 \}=1$;

\item
its joint distribution is
\begin{eqnarray*}
& & \Pr\{ W_{t_1}\in \Gamma_1, ... , W_{t_n}\in\Gamma_n \}= \\
& & \int\limits_{\Gamma_1}dx_1 ... \int\limits_{\Gamma_n}dx_n ~
g(0,x_1~;~t_1)~ g(x_1,x_2~;~t_2-t_1)~ ... ~ g(x_{n-1},x_n~;~t_n-t_{n-1}) \\
\end{eqnarray*}
for any integer $n$, any real numbers $0<t_1<...<t_n$ and arbitrary
intervals $\Gamma_1$, ..., $\Gamma_n$, where $g(x,x'~;~t)$ is the
Gaussian density
\begin{equation}                                                                              
\label{eq:Gaussian_density}
g(x,x'~;~t)=\frac{1}{\sqrt{2\pi t}}~ \exp\left[-\frac{(x-x')^2}{2t}\right]
\hskip 15mm x,x'\in\R, ~~~ t\in\R_+
\end{equation}
By definition, $g(0,x~;~t)dx$ is the probability to find the Brownian
motion in $dx$ vicinity of point $x$ at time $t$:
\begin{equation*}
\Pr\{ W_t\in (x,x+dx) \}=g(0,x~;~t)dx
\end{equation*}
\end{itemize}

The collection $W_t=(W_t^1,...,W_t^d)$ of $d$ independent one-dimensional
Brownian motions $W_t^k$ is called {\it $d$-dimensional Brownian motion}
started from the origin (in the following, we shall omit the pointing on 
the dimension). The translated stochastic process, $x+W_t$, is called 
Brownian motion started from the point $x\in\R^d$.
\end{dfntn}

Various properties of the Brownian motion and its relation to other
mathematical fields (like partially differential equations or
potential theory) are well known and can be found in
\cite{Freidlin,Borodin,Ito,Levy,Revuz,Bass,Port}.

As one can see, the Brownian motion $W_t$ is defined for the whole
space $\R^d$, without any binding to a particular domain. However,
physical processes are usually confined into a certain domain
$\Omega\subset\R^d$.  The ``presence'' of its boundary $\pa$ can be
introduced by a specific condition for a quantity we are looking
for. To illustrate this notion, let us introduce the harmonic measure
$\omega_x$ defined as the probability measure to hit different subsets
of the boundary $\pa$ for the first time.

\begin{dfntn}
Let $\Omega\subset\R^d$ be a domain with boundary $\pa$.  For any
$x\in\R^d$, a random variable $\T^x=\inf\{ t>0 ~:~ (x+W_t)\in\pa \}$
is called {\it stopping time} on the boundary $\pa$ (it gives the
first moment when the Brownian motion started from $x$ hits the
boundary).  For any subset $A$ from the Borel $\sigma$-algebra
$\B(\pa)$, one defines its harmonic measure $\omega_x\{ A \}$ (hitting
probability) as:
\begin{equation*}
\omega_x\{ A \}=\Pr\{ W_{\T^x}\in A , ~ \T^x <\infty \}
\end{equation*}
(we remind that the Borel $\sigma$-algebra $\B(\pa)$ is generated by all 
open subsets of $\pa$). 
\end{dfntn}

We gave this classical definition of the harmonic measure in order to
outline that the boundary $\pa$ is present in the problem only through
the stopping time $\T^x$. In other words, its introduction does not
change the definition of the Brownian motion itself. This feature
considerably simplifies the following analysis.

Up to this moment, we did not specify the domain $\Omega$ and its
boundary $\pa$, since the harmonic measure can be well defined for
very irregular domains \cite{Garnett,Garnett2,Makarov99}. However, the
following definitions will need some restrictions on the boundary.
Throughout this paper, we shall consider a domain $\Omega\subset\R^d$
($d\geq 2$) with {\it bounded smooth boundary} $\pa$ (twice continuous
differentiable manifold).  One the one hand, this condition can be
weakened in different ways, but it would require more sophisticated
analysis overflowing the frames of this paper (e.g., see
\cite{Grisvard,Grisvard2}).  One the other hand, our primary aim is
to describe Laplacian transport phenomena listened in
Section~\ref{sec:PTL}. Dealing with {\it physical} problems, one can
always {\it smooth} a given boundary $\pa$ whatever its original
irregularity. Indeed, the physics naturally provides a minimal cut-off
$\delta$ (e.g., mean free path of diffusing or reacting molecules)
which determines the ``admissible'' scales of the boundary. All
geometrical features of the boundary smaller than $\delta$ should be
irrelevant (otherwise, the proposed physical description would be
incorrect). Smoothing these geometrical elements, one can obtain a
boundary that may be (very) irregular on length scales larger than
$\delta$, but smooth on length scales lower than $\delta$. For a
smooth boundary $\pa$, one can introduce the harmonic measure density
$\omega_x(s)$ such that $\omega_x(s)ds$ is the probability that the
Brownian motion started at $x$ hits the boundary in $ds$ vicinity of
the boundary point $s$.

The harmonic measure, generated by the Brownian motion, gives a
general solution of the Dirichlet boundary value problem with a given
function $f$ on $\pa$:
\begin{equation}                                                                                  
\label{eq:Dirichlet}
\Delta u=0 ~~~ (x\in \Omega),  \hskip 15mm  u=f ~~~ (x\in\pa)
\end{equation}
Indeed, the harmonic measure density is equal to the normal derivative
of the Green function for the Dirichlet problem, so that one writes
the solution $u(x)$ explicitly \cite{Roach}:
\begin{equation*}
u(x)=\int\limits_\pa f(s) ~\omega_x(s)ds
\end{equation*}
or as following expectation
\begin{equation}                                                                               
\label{eq:Dirichlet_E}
u(x)=\E \bigl\{ f(W_{\T^x}) \bigr\}
\end{equation}

One can give a physical interpretation to this mathematical relation.
In order to calculate the expectation, one considers all possible
trajectories of the Brownian motion started from the point
$x\in\Omega$.  For each trajectory terminated at boundary points
$s=W_{\T^x}$, one assigns the weight $f(s)$ and then averages over all
these trajectories.  Giving this interpretation, we do not discuss the
mathematical realization of such average over all possible
trajectories. To do this operation properly, one can introduce the
Wiener measure on the space of continuous functions and then define
the corresponding functional integrals \cite{Freidlin}. It is
interesting to remark that this reasoning traced to the Feynman's
description of quantum mechanics by path integrals \cite{Feynman}.
Note also that the relation (\ref{eq:Dirichlet_E}) is the mathematical
foundation to Monte Carlo numerical tools for solving the Dirichlet
problem (\ref{eq:Dirichlet}). In fact, launching a large number of
random walkers from the point $x$, one determines, for each
trajectory, its hitting point $s$ and assigns the corresponding weight
$f(s)$. The average over all random walkers gives an approximate value
of the solution $u(x)$ at point $x$.

One can see that the Brownian motion is an efficient mathematical tool
to study Dirichlet boundary value problems. However, it becomes
useless for other types of boundary conditions like, e.g., the Neumann
condition.  The simple physical reason is the following. As we have
mentioned above, the Dirichlet boundary condition is introduced
through the stopping time $\T^x$. It means that we are interested only
in the Brownian motion $W_t$ for times $t$ between $0$ and $\T^x$.
Since the motion with $t>\T^x$ is irrelevant for this problem, one may
think that the Brownian motion is {\it absorbed} on the boundary $\pa$
at the first hit.  In other words, the Dirichlet condition corresponds
to a {\it purely absorbing} interface $\pa$. For the Neumann
condition, the situation changes drastically. The normal derivative
representing a flux leads to the notion of {\it reflection} on the
boundary: if one would like to fix the flux density across the
boundary, certain particles should be reflected. The probabilistic
description of the Neumann boundary condition necessitates thus an
introduction of the other stochastic process called {\it reflected
Brownian motion}.

\subsection{ Reflected Brownian Motion }

The fact of reflection on the boundary implies three essential
distinctions with respect to the (simple) Brownian motion:
\begin{itemize}
\item
the definition of the reflected Brownian motion will depend on the
domain $\Omega$ (as a consequence, it will be necessarily more
sophisticated than the above definition of the Brownian motion);

\item
the type and direction of each reflection should be prescribed (e.g.,
right or oblique);

\item
some restrictions on the boundary $\pa$ should be introduced, for
instance, the normal vector should be well defined at each point (as a
consequence, the boundary cannot be very irregular).
\end{itemize}
It is not thus surprising that the definition of the reflected
Brownian motion requires {\it stochastic} differential equations. We
do not intend to reproduce the whole analysis leading to the reflected
Brownian motion since one can find it in corresponding literature
(e.g., see \cite{Freidlin,Anderson76,Lions84,Saisho87}).  In the case
of smooth boundaries, the following definition is quite classical. The
situation becomes essentially more difficult when one tries to extend
it for nonsmooth domains.

\begin{dfntn}
Let $\Omega\subset\R^d$ be a domain with boundary $\pa$, and $n(s)$ is
a vector-valued function on $\pa$. For a given point $x\in \Omega$,
one considers the stochastic equation in the following form:
\begin{equation}                                                                              
\label{eq:stochastic}
d\hat{W}_t = dW_t + n(\hat{W}_t) \I_\pa(\hat{W}_t) d\l_t  \hskip 15mm
\hat{W}_0=x, ~~ \l_0=0
\end{equation}
where $W_t$ is $d$-dimensional Brownian motion and $\I_\pa$ is the
indicator of the boundary $\pa$. By a solution of this equation, we
mean a pair of almost surely continuous processes $\hat{W}_t$ and
$\l_t$, satisfying (\ref{eq:stochastic}), adapted to the underlying
family of $\sigma$-fields and satisfying, with probability $1$, the
following conditions:
\begin{itemize}
\item
$\hat{W}_t$ belongs to $\Omega\cup\pa$;

\item
$\l_t$ is a nondecreasing process which increases only for $t\in \TT$,
$\TT=\{ t>0 ~:~ \hat{W}_t\in \pa \}$ having Lebesgue measure zero
almost surely.
\end{itemize}
The process $\hat{W}_t$ is called {\it Brownian motion normally
reflected on the boundary} (or {\it reflected Brownian motion}), the
process $\l_t$ is called {\it local time on the boundary} (or {\it
local time process}).
\end{dfntn}

The following theorem ensures the existence and uniqueness of these
stochastic processes in the case of smooth boundaries.

\begin{thrm}
Let $\Omega\subset \R^d$ be a bounded domain with twice continuous
differentiable boundary $\pa$, $n(s)$ is the vector of the inward unit
normal at boundary point $s$ (orthogonal to the boundary at $s$ and
oriented towards the domain). For a given point $x\in\Omega\cup\pa$,
the stochastic equation (\ref{eq:stochastic}) possesses a unique
solution, i.e., there exist the reflected Brownian motion $\hat{W}_t$
and the local time on the boundary $\l_t$ satisfying the above
conditions, and they are unique.
\end{thrm}
\underline{Proof} can be found in \cite{Anderson76,Freidlin}.

\vskip 2mm
We should note that this theorem can be extended in different ways.
For example, one can consider the Brownian motion, reflected on the
boundary in the direction given by another vector-valued field than
the field $n(s)$ of the inward unit normals. The assumption that the
domain is bounded can be replaced by a more subtle hypothesis that
allows to extend the definition of the reflected Brownian motion for
some classes of unbounded domains. At last, one may define this motion
for a general case of second order elliptic differential operators
(with certain restrictions on their coefficients).  The interested
reader may consult the corresponding literature, e.g.,
\cite{Freidlin,Bass}.

Although the rigorous mathematical definition of stochastic
differential equations is more difficult than in the case of ordinary
differential equations, an intuitive meaning of its elements remains
qualitatively the same. For example, the stochastic equation
(\ref{eq:stochastic}) states that an infinitesimal variation
$d\hat{W}_t$ of the reflected Brownian motion $\hat{W}_t$ in the
domain $\Omega$ (bulk) is governed only by the variation $dW_t$ of the
(simple) Brownian motion $W_t$ (the second term vanishes due to the
indicator $\I_\pa$). When the motion hits the boundary, the second
term does not allow to leave the domain leading to a variation
directed along the inward unit normal $n(s)$ towards the interior of
the domain. On the other hand, each hit of the boundary increases the
local time $\l_t$. Consequently, the single stochastic equation
(\ref{eq:stochastic}) defines simultaneously two random processes,
$\hat{W}_t$ and $\l_t$, strongly dependent each of other.

As an example, one can consider one-dimensional Brownian motion
reflected at zero which can be written as mirror reflection of the
(simple) Brownian motion: $\hat{W}_t=|x+W_t|$. Applying It\^o's
formula to this function, one obtains:
\begin{equation*}
\hat{W}_t=|x+W_t|=x ~ + ~ \int\limits_0^t \sign(x+W_{t'})dW_{t'} ~+~
\frac12 \int\limits_0^t \delta(x+W_{t'})dt'
\end{equation*}
One can show that the second term is equivalent to a Brownian motion
$W'_t$, whereas the third term, denoted as $\l_t$, is a continuous,
nondecreasing random process which increases only on the set $\TT=\{
t>0 ~:~ x+W_t=0 \}$ of the Lebesgue measure zero. The previous
expression can thus be written as $\hat{W}_t=x + W'_t + \l_t$ or, in
differential form, as $d\hat{W}_t=dW'_t + d\l_t$ which is the
particular case of the stochastic equation (\ref{eq:stochastic}).  For
the local time $\l_t$, L\'evy proved the following representation
\cite{Levy,Ito}:
\begin{equation}                                                                                        
\label{eq:Levy1}
\l_t=\lim\limits_{a\to 0} \frac{1}{2a}\int\limits_0^t \I_{[0,a]}(\hat{W}_{t'})dt'
\end{equation}
This relation makes explicit the meaning of the local time $\l_t$: it
shows how ``many times'' the reflected Brownian motion passed in an
infinitesimal vicinity of zero up to the moment $t$. L\'evy also gave
another useful representation for the local time:
\begin{equation}                                                                                        
\label{eq:Levy2}
\l_t=\lim\limits_{a\to 0} a\N_t(a)   
\end{equation}
where $\N_t(a)$ is the number of passages of the reflected Brownian
motion through the interval $[0,a]$ up to the moment $t$. If one 
introduces a sequence of stopping times at points $0$ and $a$,
\begin{eqnarray*}
\tau^{(0)}_0=\inf\{ t>0 ~:~ \hat{W}_t=0 \}  && 
\tau^{(a)}_0=\inf\{ t>\tau^{(a)}_0 ~:~ \hat{W}_t=a \} \\
\tau^{(0)}_n=\inf\{ t>\tau^{(a)}_{n-1} ~:~ \hat{W}_t=0 \}  &&
\tau^{(a)}_n=\inf\{ t>\tau^{(0)}_{n-1} ~:~ \hat{W}_t=a \}
\end{eqnarray*}
the number of passages can be defined as
\begin{equation*}
\N_t(a)=\sup\{ n>0 ~:~ \tau^{(0)}_n<t \}
\end{equation*}

Note that the representations (\ref{eq:Levy1}) and (\ref{eq:Levy2})
can be extended for a general case of $d$-dimensional reflected
Brownian motion.

\section{ Partially Reflected Brownian Motion }                   
\label{sec:PRBM}

\subsection{ Definition and Certain Properties }

Bearing in mind the description of Laplacian transport phenomena, we
would like to extend the concept of the reflected Brownian motion in
order to deal with the mixed boundary condition (\ref{eq:problem_C3}).

\begin{dfntn}                                                                                       
\label{th:PRBM}
For a given domain $\Omega\subset\R^d$ with smooth bounded boundary
$\pa$, let $\hat{W}_t$ be the reflected Brownian motion started from
$x\in\Omega\cup\pa$, and $\l_t$ be the related local time process. Let
$\chi$ be a random variable, independent of $\hat{W}_t$ and $\l_t$ and
distributed according to the exponential law with a positive parameter
$\Lambda$:
\begin{equation}                                                                                    
\label{eq:expon}
\Pr\{ \chi \geq \lambda \}=\exp[-\lambda/\Lambda]   \hskip 15mm  (\lambda\geq 0)
\end{equation}
The stopping time 
\begin{equation*}
\T^x_\Lambda=\inf\{ t>0 ~:~ \l_t\geq \chi \}
\end{equation*}
gives the first moment when the local time process $\l_t$ exceeds the
random variable $\chi$. The process $\hat{W}_t$ conditioned to stop at
random moment $t=\T^x_\Lambda$ is called {\it partially reflected
Brownian motion} (PRBM).
\end{dfntn}

First of all, we stress that the partially reflected Brownian motion
is not a {\it new} stochastic process: it reproduces completely the
reflected Brownian motion $\hat{W}_t$ up to the moment $\T^x_\Lambda$.
The only difference between them resides in the fact that we are not
interested in what happens after this moment. Consequently, the
condition to stop at $t=\T^x_\Lambda$ may be thought as an absorption
on the boundary $\pa$. It explains the term ``partially reflected'':
after multiple reflections, the process will be absorbed on the
boundary (see Section~\ref{sec:approaches} for further comments).
Roughly speaking, the whole term ``partially reflected Brownian
motion'' is a shorter version of the phrase ``reflected Brownian
motion conditioned to stop at random moment $\T^x_\Lambda$''.

In the particular case $\Lambda=0$, the exponential distribution
(\ref{eq:expon}) is degenerated: $\Pr\{ \chi = 0 \}=1$ and $\Pr\{ \chi
> 0 \}=0$. Consequently, the stopping time becomes: $\T^x_0=\inf\{ t>0
~:~ \l_t>0 \}$. Since the first moment of an increase of the local
time process $\l_t$ corresponds to the first hit of the boundary
$\pa$, one obtains the stopping time of the (simple) Brownian motion:
$\T^x_0=\T^x$. One concludes that, for $\Lambda=0$, the partially
reflected Brownian motion becomes the Brownian motion conditioned to
stop at the first hit of the boundary.

To study the partially reflected Brownian motion, one can introduce a
measure quantifying absorptions on different subsets of the boundary
$\pa$.

\begin{dfntn}                                                                                       
\label{th:spread}
For any subset $A$ from the Borel $\sigma$-algebra $\B(\pa)$, one
defines its {\it spread harmonic measure} $\omega_{x,\Lambda}\{ A \}$
as:
\begin{equation*}
\omega_{x,\Lambda}\{ A \}=\Pr\{ \hat{W}_{\T^x_\Lambda}\in A , ~ \T^x_\Lambda<\infty \}
\end{equation*}
\end{dfntn}

As the harmonic measure itself, $\omega_{x,\Lambda}\{ A \}$ satisfies
the properties of a probabilistic measure, in particular,
$\omega_{x,\Lambda}\{ \pa \}=1$.  When $\Lambda$ goes to $0$, the
spread harmonic measure tends to the harmonic measure:
$\omega_{x,\Lambda}\{ A \}\to \omega_x\{ A \}$.  Since the present
definition of the PRBM requires the smoothness of the boundary, the
spread harmonic measure can be characterized by its density
$\omega_{x,\Lambda}(s)$.

\vskip 2mm
Dealing with the Brownian motion, one could formally take the starting
point $x$ on the boundary $\pa$, but it would lead to trivial results:
the stopping time $\T^x$ becomes $0$ and the harmonic measure
$\omega_x$ is degenerated to the Dirac point measure:
$\omega_x\{A\}=\I_A(x)$ (if $x\in\pa$).  In the case of the partially
reflected Brownian motion, the starting point $x$ can belong to the
domain $\Omega$ or to its boundary $\pa$: in both cases the spread
harmonic measure has nontrivial properties.

It is convenient to separate each random trajectory of the PRBM in two
parts, before and after the first hit of the boundary. The first part,
$\hat{W}_{0\leq t\leq \T^x}$, coincides with the (simple) Brownian
motion started from $x$ and conditioned to stop on the boundary, while
the second part, $\hat{W}_{\T^x\leq t\leq \T^x_\Lambda}$, coincides
with the reflected Brownian motion started on the boundary (at the
first hitting point) and conditioned to stop on the same boundary at
random moment $\T^x_\Lambda$. Since these two parts are independent,
one can write the spread harmonic measure density as
\begin{equation}                                                                                  
\label{eq:TLambda}
\omega_{x,\Lambda}(s)=\int\limits_\pa ds'\omega_x(s') ~ T_\Lambda(s',s)
\hskip 15mm  T_\Lambda(s',s)\equiv \omega_{s',\Lambda}(s)
\end{equation}

The integral kernel $T_\Lambda(s',s)$ represents the probability
density that the PRBM started from the boundary point $s'$ is stopped
(absorbed) in an infinitesimal vicinity of the boundary point $s$.
Consequently, it is sufficient to determine the probabilities of
displacements between two boundary points in order to reconstruct the
whole spread harmonic measure density.

\begin{lmm}                                                                                  
\label{th:Kakutani_Lambda}
For any subset $A$ from $\B(\pa)$ and fixed positive $\Lambda$, the
spread harmonic measure $\omega_{x,\Lambda}\{A\}$, considered as a
function of $x$, solves the mixed boundary value problem:
\begin{equation}                                                                             
\label{eq:Kakutani_Lambda}
\Delta \omega_{x,\Lambda}\{A\}=0 ~~~ (x\in\Omega), \hskip 15mm
\left[I-\Lambda \frac{\partial }{\partial n}\right]\omega_{x,\Lambda}\{A\}=\I_A(x) ~~~ (x\in\pa)
\end{equation}
\end{lmm}

\vskip 2mm
This lemma generalizes the Kakutani theorem for the harmonic measure
(when $\Lambda=0$) \cite{Kakutani44}. We do not reproduce the proof of
this lemma since it would require many technical details. It can be
also reformulated for the spread harmonic measure density:

\begin{lmm}                                                                                  
\label{th:Kakutani_Lambda2}
For any boundary point $s\in\pa$ and fixed positive $\Lambda$, the
spread harmonic measure density $\omega_{x,\Lambda}(s)$, considered as
a function of $x$, satisfies the following conditions:
\begin{equation}                                                                             
\label{eq:Kakutani_Lambda2}
\Delta \omega_{x,\Lambda}(s)=0 ~~~ (x\in\Omega), \hskip 15mm
\left[I-\Lambda \frac{\partial }{\partial n}\right]\omega_{x,\Lambda}(s)=\delta(s-x) ~~~ (x\in\pa)
\end{equation}
where $\delta(s-x)$ is the Dirac function (distribution) on the boundary.
\end{lmm}

According to this lemma, the solution of a general mixed boundary
value problem
\begin{equation*}
\Delta u=0 ~~~ (x\in\Omega), \hskip 15mm
\left[I-\Lambda \frac{\partial }{\partial n}\right]u=f ~~~ (x\in\pa)
\end{equation*}
with a given function $f$ on $\pa$ and fixed positive $\Lambda$ can be
written in two equivalent forms:
\begin{equation*}
u(x)=\int\limits_\pa f(s) ~\omega_{x,\Lambda}(s)ds = 
\E \bigl\{ f(\hat{W}_{\T^x_\Lambda}) \bigr\}
\end{equation*}
Again, one can give a physical interpretation of this relation: one
averages the function $f$ over all possible trajectories of the
partially reflected Brownian motion started from the point $x$.  Each
trajectory is weighted by $f(s)$ according to the boundary point $s$
of its final absorption.

\subsection{ Planar Surface }

We remind that the physical motivation of this work is a possibility
to describe diffusing particles near semi-permeable interfaces by the
partially reflected Brownian motion. Indeed, the mixed boundary value
problem (\ref{eq:problem_C1}--\ref{eq:problem_C3}) is an averaged
description for the concentration of particles, while the stochastic
description permits to ``follow'' the trajectory of one individual
particle. This analysis may provide a new information: typical or
average distance between the first hitting point and the final
absorption point; proportion of ``flatten'' trajectories, going near
the interface, with respect to remote trajectories, moving away from
the interface and then returning to it, etc. In this subsection, we
briefly consider the particular case of the planar surface (boundary
of a half space), when the partially reflected Brownian motion can be
constructed in a simple way, without stochastic equations.
Consequently, many related characteristics can be obtained explicitly.
In addition, this construction for the half space brings an example of
the PRBM for an unbounded domain.

Let $\Omega$ be the upper half space, $\Omega=\{ x\in\R^d ~:~ x_d>0
\}$, with smooth boundary $\pa =\{ x\in\R^d ~:~ x_d=0 \}$. Let
$\{W_t^k\}$ are $d$ independent Brownian motions started from the
origin. Then, the Brownian motion, started from a given point
$x\in\Omega$ and reflected on the boundary $\pa$, can be written in a
simple way as $(x_1+W_t^1,~ ..., ~x_{d-1}+W_t^{d-1},|x_d+W_t^d|)$. The
particular simplification is brought by the fact that reflections
happen in a single direction, being involved through the
one-dimensional reflected Brownian motion $|x_d+W_t^d|$. Without loss
of generality, we can consider the reflected Brownian motion started
from the origin ($x=0$): the translational invariance along the
hyperplane $\pa$ permits to move the starting point in $\pa$, whereas
the convolution property (\ref{eq:TLambda}) allows displacements in
orthogonal direction. The local time process $\l_t$ can be introduced
either through the stochastic equations (\ref{eq:stochastic}), or with
the help of L\'evy's formulae (\ref{eq:Levy1}) or (\ref{eq:Levy2}).

\begin{lmm}
Let $\Omega$ be the upper half space, $\Omega=\{ x\in\R^d ~:~ x_d>0
\}$.  For any positive $\Lambda$, the stopping time $\T^0_\Lambda$,
defined in~\ref{th:PRBM}, is distributed according to
\begin{equation}                                                                                 
\label{eq:temps_distrib}
\Pr\{ \T^0_\Lambda \in (t,t+dt) \}=\rho_\Lambda(t)dt ,   \hskip 20mm
\rho _\Lambda (t)=\int\limits _0^\infty \frac{z~e^{-z^2/2t}e^{-z/\Lambda }}
{\Lambda \sqrt{2\pi}~t^{3/2}}~dz 
\end{equation}
\end{lmm}
\underline{Proof}.
Since the local time $\l_t$ and the random variable $\chi$ are
independent, one can write the probability $\Pr\{ \T^0_\Lambda \in
(t,t+dt) \}$ as
\begin{equation*}
\Pr\{~\T^0_\Lambda \in (t,t+dt)~\}=\int\limits _0^\infty \Pr\biggl\{~\inf\{ \tau>0 ~:~ 
\ell_\tau=z \}\in (t,t+dt)~\biggr\} ~\Pr\{~\chi \in (z,z+dz)~\}
\end{equation*}
Then, the first factor is the well known density of the inverse local
time process \cite{Borodin},
\begin{equation*}
\Pr\biggl\{~\inf\{ \tau>0 ~:~ \ell_\tau=z \}\in (t,t+dt)~\biggr\}=dt~ 
\frac{z~e^{-z^2/2t}}{\sqrt{2\pi}~t^{3/2}}
\end{equation*}
while the second factor is given by the exponential law density
(\ref{eq:expon}) that implies (\ref{eq:temps_distrib}).~$\square$

\vskip 2mm
Note that the integral in (\ref{eq:temps_distrib}) can be represented
with the help of the Gaussian error function
\begin{equation*}
\rho _\Lambda (t)=\frac{1}{2\Lambda ^2}\left[\frac{1}{\sqrt{\pi }}~
\frac{1}{\sqrt{t/2\Lambda^2}}-\K \biggl(\sqrt{t/2\Lambda^2}\biggr) \right] ,
\hskip 15mm  \K(z)=\frac{2}{\sqrt{\pi}}\int\limits _z^\infty e^{z^2-x^2}dx
\end{equation*}
One finds the asymptotic behavior of the density $\rho_\Lambda(t)$:
\begin{equation*}
\rho_\Lambda(t)\sim \bigl(\sqrt{2\pi }\Lambda\bigr)^{-1} ~ t^{-1/2}   ~~~ (t\to 0), \hskip 20mm  
\rho_\Lambda(t)\sim \bigl(\sqrt{2\pi}/\Lambda \bigr)^{-1} ~ t^{-3/2}    ~~~ (t\to \infty)
\end{equation*}

\vskip 2mm
Once the distribution of stopping time $\T^0_\Lambda$ is determined,
one can calculate the spread harmonic measure density
$\omega_{x,\Lambda}(s)$.

\begin{lmm}
Let $\Omega$ be the upper half space, $\Omega=\{ x\in\R^d ~:~ x_d>0
\}$.  For any positive $\Lambda$, the spread harmonic measure density
$\omega_{x,\Lambda}(s)$ is
\begin{equation}                                                                                 
\label{eq:spread_planar}
\omega _{x,\Lambda}(s_1,...,s_{d-1})=\int\limits _{-\infty}^\infty
...\int\limits _{-\infty}^\infty \frac{dk_1...dk_{d-1}}{(2\pi)^{d-1}} ~
\exp\biggl[-i\sum\limits _{j=1}^{d-1}k_j(x_j-s_j)\biggr]~\frac{e^{-x_d |k|}}{1+\Lambda |k|}
\end{equation}
where $|k|=\sqrt{k_1^2+...+k_{d-1}^2}$.
\end{lmm}
\underline{Proof}.
First, the probability kernel $T_\Lambda(s,s')$, defined for two
boundary points $s,s'\in\pa$, is translationally invariant in the
hyperplane $\pa$, $T_\Lambda(s,s')=t_\Lambda(s-s')$, where
\begin{equation*}
\begin{split}
t_\Lambda (s_1,...,s_{d-1})&ds_1...ds_{d-1}\equiv \\
& \Pr\biggl\{~W^1_{\T^0_\Lambda} \in (s_1,s_1+ds_1)~, 
~... ~,~ W^{d-1}_{\T^0_\Lambda} \in (s_{d-1},s_{d-1}+ds_{d-1})~\biggr\} \\
\end{split}
\end{equation*}
The stopping time $\T^0_\Lambda$ is related to the orthogonal motion
and, consequently, independent of lateral motions $W_t^1$, ...,
$W_t^{d-1}$.  Therefore, the above probability can be written as
\begin{equation*}
\begin{split}
t_\Lambda &(s_1,...,s_{d-1})ds_1...ds_{d-1}= \\
& \int\limits _0^\infty \Pr\biggl\{~
W^1_t\in (s_1,s_1+ds_1), ~ ... ~,~ W^{d-1}_t\in (s_{d-1},s_{d-1}+ds_{d-1})~\biggr\}~\rho_\Lambda(t)dt \\
\end{split}
\end{equation*}
Since the lateral motions are independent between themselves, the
first factor is equal to the product of Gaussian densities
(\ref{eq:Gaussian_density}):
\begin{equation*}
t_\Lambda(s_1,...,s_{d-1})=\int\limits _0^\infty dt~ \rho_\Lambda (t) ~
\prod\limits_{j=1}^{d-1} \frac{e^{-s^2_j/2t}}{\sqrt{2\pi t}}
\end{equation*}
Using the integral representation (\ref{eq:temps_distrib}), one finds
\begin{equation*}
t_\Lambda (s_1,...,s_{d-1})=\frac{\Gamma (d/2)}{\pi ^{d/2}~\Lambda }\int\limits _0^\infty
dz ~\frac{z~e^{-z/\Lambda }}{\bigl[s_1^2+...+s_{d-1}^2+z^2\bigr]^{d/2}}
\end{equation*}
where $\Gamma(z)$ stands for Euler gamma function (see
\cite{Grebenkov} for details).

Substituting the well known harmonic measure density $\omega_x(s)$ for
the upper half space (generalized Cauchy distribution),
\begin{equation*}
\omega _x(s_1,...,s_{d-1})=\frac{\Gamma (d/2)}{\pi ^{d/2}}~
\frac{x_d}{\bigl[(x_1-s_1)^2+...+(x_{d-1}-s_{d-1})^2+(x_d)^2\bigr]^{d/2}}
\end{equation*}
into convolution (\ref{eq:TLambda}), one finally obtains the
expression (\ref{eq:spread_planar}) for the spread harmonic measure
density.~$\square$

\vskip 2mm
One can easily verify that the spread harmonic measure density
$\omega_{x,\Lambda}(s)$ and the probability kernel $T_\Lambda(s,s')$
satisfy the following conditions:
\begin{enumerate}
\item
Normalization condition:
\begin{equation*}
\int\limits_{\pa} ds ~\omega_{x,\Lambda}(s) = 1   \hskip 20mm  
\int\limits_{\pa} ds' ~T_\Lambda(s,s') = 1
\end{equation*}

\item
Dirichlet limit ($\Lambda\to 0$):
\begin{equation*}
\omega_{x,\Lambda}(s)\longrightarrow \omega_x(s)   \hskip 20mm
T_\Lambda(s,s')\longrightarrow \delta(s-s')
\end{equation*}

\item
Translational invariance:
\begin{equation*}
\omega_{x,\Lambda}(s)=\omega_{x-s,\Lambda}(0)   \hskip 20mm
T_\Lambda(s,s')=T_\Lambda(s-s',0)\equiv t_\Lambda(s-s')
\end{equation*}

\end{enumerate}
One can also deduce the asymptotic behavior of the function
$t_\Lambda(s)$ as $|s|\to 0$ or $|s|\to \infty$. For this purpose, it
is convenient to define the new function $\eta_d(z)$ by relation
\begin{equation*}
t_\Lambda(s)=\eta _d\bigl(|s|/\Lambda\bigr) ~ \omega_{(0,...,0,\Lambda)}(s)
\end{equation*}
where the second factor is the harmonic measure density for the
Brownian motion started from the point
$(\underbrace{0,...,0}_{d-1},\Lambda)$.  Using the explicit formulae
for $t_\Lambda(s)$ and $\omega_{(0,...,0,\Lambda)}(s)$, one obtains:
\begin{equation}
\eta_d(z)=\bigl(1+z^2\bigr)^{d/2}\int\limits _0^\infty \frac{t~ e^{-t}~ dt}{(t^2+z^2)^{d/2}} 
\end{equation}
Its asymptotic behavior for $z$ going to infinity is
\begin{equation}                                                                          
\label{eq:etad_inf}
\eta _d(z)=1-\frac{5d}{2}~ z^{-2} + O(z^{-4})
\end{equation}
%%  d(109d+238)/8 z^{-4}
whereas for $z$ going to $0$, one has
\begin{equation*}
\eta _d(z)\sim \frac{\Gamma(d/2)}{\pi ^{d/2}}~z^{2-d}  ~~~ (d>2), 
\hskip 15mm  \eta _d(z)\sim \frac{1}{\pi} ~\ln z  ~~~ (d=2)
\end{equation*}

These relations can be used for qualitative study of the partially
reflected Brownian motion. For instance, one identifies the parameter
$\Lambda$ as a characteristic length scale of the problem: the
magnitude of any distance (e.g., $|s|$) has to be compared with
$\Lambda$.  Interestingly, the asymptotic behavior (\ref{eq:etad_inf})
for large $z$ means that the function $t_\Lambda(s)$ is close to the
harmonic measure density $\omega_{(0,...,0,\Lambda)}(s)$. Roughly
speaking, for large $|s|/\Lambda$, the partially reflected Brownian
motion started from the origin is qualitatively equivalent to the
(simple) Brownian motion started from the point $(0,...,0,\Lambda)$.
In other words, the partial reflections on the boundary lead to a {\it
spreading} of the harmonic measure with characteristic scale $\Lambda$
(see also relation (\ref{eq:TLambda})). The explicit analytical
results can also be derived in presence of an absorbing barrier at a
given height \cite{Sapoval05}.

\vskip 2mm
The knowledge of the probability kernel $T_\Lambda(s,s')$ brings an
important information about the partially reflected Brownian motion in
the (upper) half space. As an example, we calculate the probability
$\Prob_\Lambda(r)$ that the PRBM started from the origin is finally
absorbed on the disk $B^{d-1}_r=\{~ (x_1,...,x_d)\in\R^d ~:~
x_1^2+...+x_{d-1}^2\leq r^2 ,~ x_d=0 ~\}$ of radius $r$ centered at
the origin:
\begin{equation}                                                                                   
\label{eq:PBLambda}
\Prob_\Lambda(r)=\int\limits_{B^{d-1}_r} ds~ t_\Lambda(s)=
\frac{2~\Gamma(\frac{d}{2})}{\Gamma (\frac{d-1}{2})\sqrt{\pi}}\int\limits _0^\infty 
te^{-t}dt \int\limits _0^{r/\Lambda} \frac{ x^{d-2}~dx}{[x^2+t^2]^{d/2}} 
\end{equation}

This probability shows how far the partially reflected Brownian motion
can go away after the first hit of the boundary. One sees that this
function depends only on the ratio $r/\Lambda$, going to $0$ for small
radii and to $1$ for large radii. Again, the parameter $\Lambda$ is
the characteristic length scale of the problem. In two-dimensional
case, $\Prob_\Lambda(\Lambda/2)$ is the probability that the PRBM is
absorbed on the linear segment of length $\Lambda$, centered at the
origin (the first hitting point). The numerical calculation of the
integral in (\ref{eq:PBLambda}) gives $\Prob_\Lambda(\Lambda/2)\simeq
0.4521$, i.e., about half of the particles is absorbed on this
region. In other words, the length of the characteristic absorption
region (where half of the particles is absorbed) is approximately
equal to $\Lambda$. It has been shown recently that this result is
qualitatively valid for a large class of irregular boundaries
\cite{Grebenkov06c}.  Roughly speaking, if the one-dimensional
boundary (curve) has no deep pores (fjords) and its perimeter is large
with respect to the scale $\Lambda$, then the curvilinear interval of
length $\Lambda$, centered on the first hitting point, absorbs
approximately half of the diffusing particles.

\begin{figure}
\includegraphics[width=138mm]{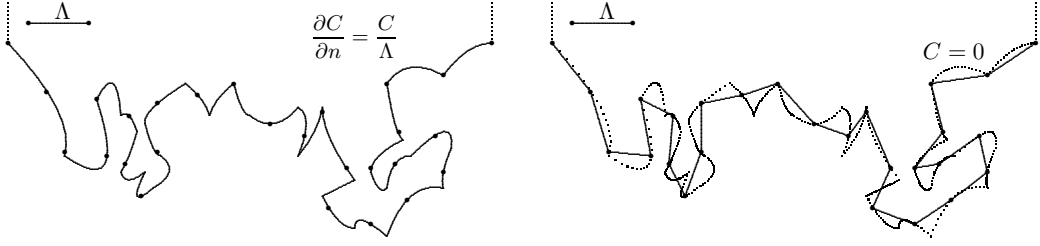}
\caption{ 
Land Surveyor Approximation: the total flux across the boundary can be
approximately calculated when the mixed boundary condition $\partial
C/\partial n =C/\Lambda$ on a given irregular curve (on the left) is
replaced by the Dirichlet condition $C=0$ on the coarse-grained
boundary (on the right). The last one is obtained by replacing
curvilinear intervals of length $\Lambda$ by corresponding linear
chords. }
\label{fig:LSA}
\end{figure}

This result can be considered as a first mathematical justification of
the {\it Land Surveyor Approximation} (LSA) developed by Sapoval
\cite{Sapoval94}.  According to this approximation, a given
one-dimensional interface (curve) can be coarse-grained with physical
scale $\Lambda$ in order to replace the mixed boundary condition
$[I-\Lambda \partial /\partial n]C=0$ by the Dirichlet condition $C=0$
(see Fig.~\ref{fig:LSA}). Some heuristic physical arguments allowed to
state that the total flux across the irregular semi-permeable
interface was approximately equal to the total flux across this
coarse-grained boundary with Dirichlet condition. This statement
provided a simple but powerful tool to investigate Laplacian transport
phenomena. The land surveyor approximation had been checked
numerically \cite{Filoche97,Sapoval99}, but not mathematically.  The
study of the partially reflected Brownian motion brings its
justification and further understanding. Actually, the coarse-graining
procedure generates the regions of length $\Lambda$, where about half
of the particles is absorbed.  The LSA is based on two simplifications
which can be clearly explained in terms of diffusing particles:
\begin{enumerate}
\item
The Dirichlet boundary condition on the coarse-grained boundary means
that all particles arrived to the characteristic absorption region are
absorbed. This approximation does not take into account half of the
particles which escaped this region.

\item
The linear chords, generated by coarse-graining, are deterministic
regions. This approximation neglects the fact that the characteristic
absorption regions should be centered at the {\it random} position of
the first hit of the boundary.

\end{enumerate}
Although these simplifications seem to be rough, the numerical
simulations show that the LSA reproduces the transport properties with
good accuracy. However, this approximation has no any kind of small
parameter which would allow to control its applicability. More
accurate theoretical approaches will be discussed in 
Section~\ref{sec:approaches}.

One can go further by extending the Land Surveyor Approximation to the
three-dimensional case, which still remains poorly understood. Indeed,
the numerical calculation in 3D leads to $\Prob_\Lambda(\Lambda)\simeq
0.4611$, i.e., about half of the particles is absorbed on the disk of
radius $\Lambda$ centered at the first hitting point. Consequently, if
one finds a convenient cover of a semi-permeable irregular interface
by disk-like sets of characteristic radius $\Lambda$, the LSA may be
still valid, i.e., the total flux across a given interface would be
approximated by the total flux across the perfectly absorbing
coarse-grained interface (with Dirichlet condition). An accurate
mathematical formulation of this extension and its numerical
verification present open interesting problems.

\subsection{ Relation to the Dirichlet-to-Neumann Operator }

The construction of the partially reflected Brownian motion for a
given domain $\Omega$ requires the resolution of the stochastic
differential equation (\ref{eq:stochastic}), a quite difficult
problem. Fortunately, many characteristics of this process, e.g., the
spread harmonic measure density $\omega_{x,\Lambda}(s)$, can be
obtained in another way.  This subsection is devoted to the
Dirichlet-to-Neumann operator and its relation to the partially
reflected Brownian motion.

\begin{dfntn}
For a given domain $\Omega\subset\R^d$ ($d\geq 2$) with smooth bounded
boundary $\pa$, let $u:\Omega\cup\pa\to\R$ be a harmonic function with
Dirichlet condition $u=f$, a function $f$ being from the Sobolev space
$H^1(\pa)$ (in other words, $u$ is the solution of the boundary value
problem (\ref{eq:Dirichlet})). Applying the normal derivative to $u$,
one obtains a new function $g=\partial u/\partial n$ belonging to the
space $L^2(\pa)$ of measurable and square integrable functions. Then
the operator $\M$, acting from $H^1(\pa)$ to $L^2(\pa)$, which
associates the new function $g$ with a given $f$, is called {\it
Dirichlet-to-Neumann operator}.
\end{dfntn}

It is known that the Dirichlet-to-Neumann operator $\M$ is
self-adjoint pseudodifferential operator of the first order, with
discrete positive spectrum $\{ \mu_\alpha \}$ and smooth
eigenfunctions forming a complete basis in $L^2(\pa)$
\cite{Agranovich97,Birman,Hormander,Egorov,Jacob,Schulze,Taylor,Taylor2,Grisvard,Grisvard2}.
One can also define its resolvent operator $T_\Lambda=[I+\Lambda
\M]^{-1}$, called {\it spreading operator}. This is an analytic
operator function in the whole complex plane, except a denumerable set
of points, $\C\backslash \{ -\mu_\alpha^{-1} \}$. In particular,
$T_\Lambda$ is well defined for any positive $\Lambda$.

\begin{lmm}
For any strictly positive $\Lambda$, the spreading operator
$T_\Lambda$ acts from $L^2(\pa)$ to $L^2(\pa)$ as a compact integral
operator,
\begin{equation*}
[T_\Lambda f](s)=\int\limits_\pa ds' ~ f(s') ~ T_\Lambda(s',s)
\end{equation*}
where the kernel $T_\Lambda(s,s')$ is given by (\ref{eq:TLambda}).
\end{lmm}
\underline{Proof}. 
The probability kernel $T_\Lambda(s,s')$ is a positive function
satisfying the normalization:
\begin{equation*}
\int\limits_\pa T_\Lambda(s,s')ds' = 1
\end{equation*}
since the partially reflected Brownian motion is conditioned to be
finally absorbed on the boundary. Therefore, one obtains:
\begin{equation*}
\int\limits_\pa \int\limits_\pa ds ~ds' ~|T_\Lambda(s,s')|^2 = S_{tot} <\infty
\end{equation*}
where $S_{tot}$ is the total surface area of the boundary $\pa$.  The
integral operator $T_\Lambda$ defined by the kernel $T_\Lambda(s,s')$
is a Hilbert-Schmidt operator and, consequently, a compact operator.

The boundary condition in lemma~\ref{th:Kakutani_Lambda2} can be
written with the help of the Dirichlet-to-Neumann operator:
\begin{equation*}
\bigl[I+\Lambda \M]T_\Lambda(s,s')=\delta(s-s')
\end{equation*}
that implies that the integral operator $T_\Lambda$ coincides with the
resolvent $[I+\Lambda \M]^{-1}$ of the Dirichlet-to-Neumann operator
$\M$.~$\square$

\vskip 2mm
This simple lemma creates a ``bridge'' between the partially reflected
Brownian motion and the Dirichlet-to-Neumann operator. In particular,
the relation (\ref{eq:TLambda}) for the spread harmonic measure
density $\omega_{x,\Lambda}(s)$ can now be understood as application
of the spreading operator $T_\Lambda$ to the harmonic measure density
$\omega_x(s)$. Consequently, once the Dirichlet-to-Neumann operator is
constructed for a given domain, one can calculate the density
$\omega_{x,\Lambda}(s)$ without solving the stochastic differential
equations (\ref{eq:stochastic}).

The self-adjointness of the Dirichlet-to-Neumann operator allows one
to apply efficient tools of the spectral theory. For example, one can
rewrite the relation (\ref{eq:TLambda}) as spectral decomposition of
the harmonic measure density on eigenfunctions $\V_\alpha$ of the
operator $\M$:
\begin{equation}                                                                                    
\label{eq:omega_M}
\omega_{x,\Lambda}(s)=\sum\limits_\alpha \frac{\bigl(\omega_x \cdot \V_\alpha^* \bigr)_{L^2}}
{1+\Lambda \mu_\alpha}~ V_\alpha(s)  \hskip 15mm \bigl(\omega_x \cdot \V^*_\alpha \bigr)_{L^2}=
\int\limits_\pa \omega_x(s')~ \V^*_\alpha(s')~ ds'
\end{equation}
where $(~\cdot~)_{L^2}$ denotes the scalar product in $L^2(\pa)$
space.  The advantage of this relation is an explicit dependence on
the physical parameter $\Lambda$.

\subsection{ Examples }

In order to illustrate the underlying concepts, we consider several
examples.

\vskip 2mm
\subsubsection{ Two-Dimensional Disk }                                                                  
\label{sec:disk}

We are going to study the partially reflecting Brownian motion in a
unit disk, $\Omega =\{~ x\in \R^2 ~:~ |x|<1 ~\}$ (its boundary is a
unit circle, $\pa =\{~ x\in \R^2 ~:~ |x|=1 ~\}$).

In this case, the harmonic measure density $\omega_x(s)\equiv\omega
(r,\theta)$ is a function of two real variables: the distance $0\leq
r<1$ between the starting point $x\in \Omega$ and the origin, and the
angle $0\leq \theta <2\pi $ between directions onto points $x$ and
$s\in\pa$ from the origin. The harmonic measure density is known as
Poisson kernel:
\begin{equation}                                                                                    
\label{eq:harm_disk}
\omega (r,\theta )=\frac{1-r^2}{2\pi (1-2r\cos \theta +r^2)}
\end{equation}

The rotational invariance of the domain $\Omega $ implies that the
eigenbasis of the Dirichlet-to-Neumann operator $\M$ is the Fourier
basis,
\begin{equation*}
\V_\alpha (\theta )=\frac{e^{i\alpha \theta}}{\sqrt{2\pi}} \hskip 10mm (\alpha \in \Z)
\end{equation*}
Taking Fourier harmonic as boundary condition, $u(r=1,\theta) =
e^{i\alpha \theta}$, one finds a regular solution of the corresponding
Dirichlet problem: $u(r,\theta)=r^{|\alpha|} e^{i\alpha \theta}$.
Since the normal derivative coincides with the radius derivative, one
obtains the eigenvalues of the Dirichlet-to-Neumann operator:
\begin{equation*}
\mu _\alpha =|\alpha |   \hskip 10mm (\alpha \in \Z)
\end{equation*}
These eigenvalues are doubly degenerated (expect $\mu_0=0$).

The spread harmonic measure density is given by relation
(\ref{eq:omega_M}):
\begin{equation*}
\omega _{x,\Lambda}(s)\equiv \omega _\Lambda (r,\theta)=\frac{1}{2\pi}
\sum\limits _{\alpha =-\infty}^{\infty}
\frac{r^{|\alpha|} ~ e^{i\alpha \theta } }{1+\Lambda |\alpha|}
\end{equation*}
(the scalar product of the harmonic measure density $\omega_x(s)$ and
eigenfunctions $\V_\alpha^*(\theta)$ is shown to be equal to
$r^{|\alpha|}$, with $r=|x|$). In the case $\Lambda =0$, one retrieves
the Poisson representation for the harmonic measure density
(\ref{eq:harm_disk}) just as required.  The kernel of the resolvent
operator $T_\Lambda $ is
\begin{equation*}
T_\Lambda (\theta,\theta')=\frac{1}{2\pi }\sum\limits _{\alpha =-\infty}^\infty 
\frac{e^{i\alpha (\theta -\theta')}}{1+\Lambda |\alpha|}
\end{equation*}

For the exterior problem, when $\Omega=\{ x\in\R^2 ~:~ |x|>1\}$, one
obtains exactly the same results.

\vskip 2mm
\subsubsection{ Three-Dimensional Ball }

The similar arguments can be applied for higher dimensions. For
example, in the three-dimensional case, one considers the unit ball
$\Omega =\{~ x\in \R^3 ~:~ |x|<1 ~\}$. The harmonic measure density is
known to be
\begin{equation*}
\omega_x(s)\equiv \omega(r,\theta )\equiv \frac{1-r^2}{4\pi \bigl(1-2r\cos \theta +r^2\bigr)^{3/2}} 
\hskip 10mm (s\in \pa)
\end{equation*}
where $r=|x|<1$ is the distance between the starting point $x\in
\Omega $ and the origin, and $\theta $ is the angle between directions
onto points $x$ and $s\in\pa$ from the origin. This function can be
expanded on the basis of spherical harmonics as
\begin{equation*}
\omega_x(s)=\sum\limits _{l=0}^\infty \sum\limits_{m=-l}^l  r^l ~Y_{l,m}(s) ~Y_{l,m}(x/r)
\end{equation*}

The rotational symmetry of the problem implies that the eigenbasis of
the Dirichlet-to-Neumann operator is formed by spherical harmonics
$Y_{l,m}$.  A regular solution of the Dirichlet problem
(\ref{eq:Dirichlet}) in the unit ball can be written in spherical
coordinates $r$, $\theta$ and $\varphi$ as
\begin{equation*}
u(r,\theta,\varphi)=\sum\limits_{l=0}^\infty \sum\limits_{m=-l}^l f_{l,m}~ r^l ~Y_{l,m}(\theta,\varphi)
\end{equation*}
where $f_{l,m}$ are coefficients of the expansion of a given boundary
function $f$ (Dirichlet condition) on the complete basis of spherical
harmonics. Since the normal derivative coincides with the radius
derivative, one obtains
\begin{equation*}
[\M f](\theta,\varphi)=\left(\frac{\partial u}{\partial n}\right)_\pa =
\left(\frac{\partial u}{\partial r}\right)_{r=1}=
\sum\limits_{l=0}^\infty \sum\limits_{m=-l}^l f_{l,m}~ l ~ Y_{l,m}(\theta,\varphi)
\end{equation*}
i.e., the eigenvalues of the Dirichlet-to-Neumann operator $\M$ are
\begin{equation}                                                                                         
\label{eq:mu_l}
\mu_l=l   \hskip 10mm  (l\in\{0,1,2,...\}) 
\end{equation}
Note that the $l$-th eigenvalue is degenerated $n_l=(2l+1)$ times.

Interestingly%
\footnote{ 
The author thanks Dr. S.~Shadchin for valuable discussions on this
relation.},
%footnote
the Dirichlet-to-Neumann operator $\M$ for the unit ball $\Omega=\{
x\in \R^3 ~:~ |x|<1 \}$ coincides with an operator introduced by Dirac
in quantum mechanics \cite{Dirac}.  It is known that the hydrogen atom
is described by three quantum numbers: the main quantum number $n$,
the orbital quantum number $l$ and magnetic quantum number $m$.  Two
last numbers are associated with indices of spherical harmonics. Thus,
the Dirichlet-to-Neumann operator in the ball is apparently the
orbital quantum number operator for the hydrogen atom. In particular,
the degeneracy of eigenvalues of this operator can be understood from
the point of view of spin degeneracy.

The spread harmonic measure density $\omega_{x,\Lambda}(s)$ and the
spreading operator kernel $T_\Lambda(s,s')$ can be written explicitly
as spectral decompositions on the eigenbasis of the
Dirichlet-to-Neumann operator $\M$ as in the two-dimensional case.

\vskip 2mm
The eigenvalues of the Dirichlet-to-Neumann operator for
$d$-dimensional unit ball are still given by (\ref{eq:mu_l}) with
degeneracy
\begin{equation*}
n_l=\frac{(2l+d-2)}{(d-2)}~\frac{(l+d-3)!}{(d-3)! ~l!}
\end{equation*}

\vskip 2mm
The exterior problem for $\Omega =\{ x\in \R^3 ~:~ |x|>1 \}$ can be
considered in the same manner. The harmonic measure density is
$\tilde{\omega}_x(s)=(1/r)~\omega(1/r,\theta)$, with $r=|x|>1$. Using
the same expansion on spherical harmonics, one obtains $\mu_l=l+1$
with $l\in \{0,1,2,...\}$.  In particular, the lowest eigenvalue
$\mu_0=1$ is strictly positive. This difference with respect to the
spectrum $\mu_l$ for the interior problem has a simple probabilistic
origin: the Brownian motion in three dimensions is transient, i.e.,
there is a positive probability (equal to $1-1/r$) to never return to
the ball. Another explication follows from the theory of boundary
value problems for elliptic differential operators: the exterior
Neumann problem has a unique solution, while the solution of the
interior Neumann problem is defined up to a constant. Consequently,
the Dirichlet-to-Neumann $\M$ operator should be invertible for the
exterior problem that implies a simple condition for its eigenvalues:
$\mu_\alpha \ne 0$. On the contrary, $\M$ is not invertible for the
interior problem providing the condition that at least one eigenvalue
is zero.

\section{ Stochastic Approaches to Laplacian Transport Phenomena }                               
\label{sec:approaches}

In this section, we return to the Laplacian transport phenomena,
discussed at the beginning. First, we are going to introduce the
notion of source, diffusing particles started from. The definition of
the partially reflected Brownian motion will require only a minor
modification. After that, a recently developed continuous approach
will be presented with a special emphasis on its physical
significance.  Finally, we shall mention two other physical
descriptions which can now be considered as useful approximations to
the continuous approach.

\subsection{ Notion of Source }

The description of the partially reflected Brownian motion given in
the previous section does not involve a source, an important element
for Laplacian transport phenomena. In this subsection, we are going to
discuss the extension of previous definitions in order to introduce
the source. As one will see, a minor modification will be sufficient.

Throughout this subsection, we consider a bounded domain $\Omega$ with
twice continuous differentiable boundary composed of two disjoint
parts, $\pa$ and $\pa_0$, referred to as {\it working interface} and
{\it source} respectively%
\footnote{
Previously, the whole boundary of the domain had been considered as
the working interface. For this reason, we preserve the same notation
$\pa$ for this object and hope that it will not lead to ambiguities.}.
%footnote
In practice, the working interface and the source are well separated
in space, therefore one may think about a circular ring as generic
domain.

\vskip 2mm
As previously, one considers the reflected Brownian motion
$\hat{W}_t$, started from any point $x\in\Omega\cup\pa\cup\pa_0$ and
reflected on the whole boundary $\pa\cup \pa_0$, the corresponding
local time process $\l_t$, and the stopping time $\T^x_\Lambda$
defined in~\ref{th:PRBM}. Let us introduce a new stopping time $\tau$
as the first moment when the process $\hat{W}_t$ hits the source
$\pa_0$:
\begin{equation*}
\tau=\inf\{ t>0 ~:~ \hat{W}_t \in \pa_0 \}
\end{equation*}
Then, the spread harmonic measure can be introduced for any subset $A$
from Borel $\sigma$-algebra $\B(\pa)$ (defined on the working
interface alone!) as
\begin{equation}                                                                          
\label{eq:omega_Lambda_source}
\omega_{x,\Lambda}\{A\}=\Pr\{ \hat{W}_{\T^x_\Lambda}\in A , ~ \T^x_\Lambda < \tau <\infty \}
\end{equation}

We outline two distinctions with respect to the previous
definition~\ref{th:spread}:
\begin{itemize}
\item
The measure is considered on Borel subsets of the working interface
$\pa$ only, whereas the reflected Brownian motion $\hat{W}_t$ and the
local time process $\l_t$ are defined on the whole boundary $\pa \cup
\pa_0$.

\item
There is a supplementary condition $\T^x_\Lambda <\tau$ providing that
the process $\hat{W}_t$ should be stopped (absorbed) on the working
interface before hitting the source.
\end{itemize}

Note that $1-\omega_{x,\Lambda}\{\pa\}$ is the probability that the
process $\hat{W}_t$ started from a given point $x\in\Omega$ hits the
source $\pa_0$ before its final absorption on the working interface
$\pa$.

One can easily extend the lemma~\ref{th:Kakutani_Lambda} to this
spread harmonic measure:

\begin{lmm}                                                                            
\label{th:Kakutani_Lambda3}
For any subset $A$ from $\B(\pa)$ and fixed positive $\Lambda$, the
spread harmonic measure $\omega_{x,\Lambda}\{A\}$, considered as a
function of $x$, solves the boundary value problem:
\begin{equation}
\Delta \omega_{x,\Lambda}\{A\}=0 ~~~ (x\in\Omega),  \hskip 15mm  
\begin{array}{l} \displaystyle
\biggl[I-\Lambda\frac{\partial}{\partial n}\biggr]\omega_{x,\Lambda}\{A\}=\I_A(x) ~~~ (x\in\pa) \\
 ~~~~ \hskip 15mm \omega_{x,\Lambda}\{A\}=0  \hskip 10mm (x\in\pa_0)
\end{array}
\end{equation}
\end{lmm}
\underline{Proof} is similar to that of the lemma~\ref{th:Kakutani_Lambda}.
The last condition holds since $x\in\pa_0$ implies $\tau=0$.~$\square$

\begin{crllr}                                                                          
\label{th:problem_solution}
Function $C_\Lambda(x)=C_0(1-\omega_{x,\Lambda}\{\pa\})$ solves the
boundary value problem (\ref{eq:problem_C1}--\ref{eq:problem_C3}).
\end{crllr}
\underline{Proof} is a direct verification.~$\square$

Consequently, a simple introduction of the source allows one to apply
the previous description of the partially reflected Brownian motion to
study Laplacian transport phenomena. Due to reversibility of the
Brownian motion, one may think that $(1-\omega_{x,\Lambda}\{\pa\})dx$
gives also the probability to find the partially reflected Brownian
motion, started from the absorbing source, in $dx$ vicinity of the
point $x\in\Omega$, under partially absorbing condition on the working
interface $\pa$.  Note that such way of reasoning, being intuitive and
useful, is quite formal. In particular, the (simple) Brownian motion
started from the source returns to it infinitely many times with
probability $1$. If one really needs to define such a process, the
starting point should be taken slightly above the source.

\vskip 2mm
Since the boundary $\pa$ is supposed to be smooth, one can introduce
the spread harmonic measure density $\omega_{x,\Lambda}(s)$.  In turn,
the kernel of the spreading operator is defined as previously,
$T_\Lambda(s,s')\equiv \omega_{s,\Lambda}(s')$. In particular, one
retrieves the relation (\ref{eq:TLambda}):
\begin{equation*}
\omega_{x,\Lambda}(s)=\int\limits_\pa ds'\omega_{x,0}(s') ~T_\Lambda(s',s)
\end{equation*}
where the harmonic measure density $\omega_{x,0}(s)$ is defined by
relation (\ref{eq:omega_Lambda_source}) with $\Lambda=0$.

\vskip 2mm
The definition of the Dirichlet-to-Neumann operator can also be
extended to domains with a source. For a given function $f\in
H^1(\pa)$, one solves the Dirichlet problem in the domain $\Omega$:
\begin{equation*}
\Delta u=0 ~~~ (x\in\Omega),  \hskip 15mm  
\begin{array}{l} \displaystyle  u=f ~~~ (x\in\pa)  \\ 
                                u=0 ~~~ (x\in\pa_0) \end{array}
\end{equation*}
(in principle, one could consider another function on the source). For
a given function $f$ on $\pa$, the Dirichlet-to-Neumann operator $\M$,
acting from $H^1(\pa)$ to $L^2(\pa)$, associates the new function
$g=\partial u/\partial n$ on the working interface $\pa$. One can
prove general properties of this operator and its relation to the
partially reflected Brownian motion in a straight way. In particular,
the spreading operator $T_\Lambda$, defined by its kernel
$T_\Lambda(s,s')$, coincides with the resolvent operator
$[I+\Lambda\M]^{-1}$. However, some normalization properties may be
changed. In particular, for the probability kernel $T_\Lambda(s,s')$,
one has
\begin{equation*}
\int\limits_\pa T_\Lambda(s,s')ds' < 1
\end{equation*}
since the PRBM started from the working interface $\pa$ can now be
absorbed on the source.

\subsection{ Continuous Approach }

The stochastic treatment by means of the partially reflected Brownian
motion brings the solution to the problem
(\ref{eq:problem_C1}--\ref{eq:problem_C3}) describing Laplacian
transport phenomena: $C_\Lambda(x)=C_0(1-\omega_{x,\Lambda}\{\pa\})$
(see corollary \ref{th:problem_solution}).  One can go further using
the close relation to the Dirichlet-to-Neumann operator
\cite{Grebenkov06b,Grebenkov}.  According to the
lemma~\ref{th:Kakutani_Lambda3}, the density
$\omega_{x,\Lambda}\{\pa\}$, considered as a function of $x$, solves
the boundary value problem:
\begin{equation*}
\Delta \omega_{x,\Lambda}\{\pa\}=0 ~~~ (x\in\Omega),  \hskip 15mm  
\begin{array}{l} \displaystyle 
\left[I-\Lambda ~\frac{\partial }{\partial n}\right]\omega_{x,\Lambda}\{\pa\}=1 
~~~ (x\in\pa)  \\ 
~ \hskip 19.5mm \omega_{x,\Lambda}\{\pa\}=0 ~~~ (x\in\pa_0) \end{array}
\end{equation*}
The restriction of the function $\omega_{x,\Lambda}\{\pa\}$ on $\pa$
can be written with the help of the Dirichlet-to-Neumann operator $\M$
as
\begin{equation*}
\omega_{s,\Lambda}\{\pa\}=\bigl[(I+\Lambda \M)^{-1} 1\bigr](s)=[T_\Lambda 1](s)
\end{equation*}
where $1$ stands for a constant (unit) function on the working
interface.

One defines then the flux density $\phi_\Lambda(s)$ across the working
interface $\pa$:
\begin{equation*}
\phi_\Lambda(s)=D~ \frac{\partial C_\Lambda}{\partial n}(s)=-DC_0 ~ 
\frac{\partial \omega_{x,\Lambda}\{\pa\}}{\partial n}(s)
\end{equation*}
Since the normal derivative of a harmonic function can be represented
as the application of the Dirichlet-to-Neumann operator to the
restriction of this function on the boundary, one writes
\begin{equation*}
\phi_\Lambda(s)=DC_0 ~ [\M \omega_{s,\Lambda}\{\pa\}](s) = DC_0 ~ [\M T_\Lambda 1](s)
\end{equation*}
(the sign is changed due to particular orientation of the normal
derivative).  Taking $\Lambda=0$, one finds $\phi_0(s)=DC_0~[\M 1](s)$
and finally
\begin{equation*}
\phi_\Lambda(s)=[T_\Lambda \phi_0](s)
\end{equation*}

The transport properties of the working interface can be characterized
by a physical quantity called {\it spectroscopic impedance}. We remind
that the impedance of an electric scheme is defined as the tension
applied between two external poles, divided by the total electric
current passing through.  The formal analogy between the electric
problem and the diffusive transport, discussed in
Section~\ref{sec:PTL}, leads to a natural definition of the impedance
in our case as the concentration $C_0$ on the source $\pa_0$ divided
by the total flux across the working interface $\pa$:
\begin{equation*}
Z_{cell}(\Lambda)=\frac{C_0}{\int\limits_\pa ds ~\phi_\Lambda(s)}
\end{equation*}

Taking $\Lambda=0$, one deals with a purely absorbing interface $\pa$:
any particle arrived to $\pa$ is immediately absorbed (without
reflections).  In other words, such interface has no resistance for
passage across it.  Consequently, the impedance $Z_{cell}(0)$
represents the ``access resistance'' by the bulk: the possibility
that the Brownian motion can return to the source without hitting the
working interface. The resistance of the working interface alone can
thus be characterized by the difference between $Z_{cell}(\Lambda)$
and $Z_{cell}(0)$, called {\it spectroscopic impedance}:
\begin{equation*}
Z_{sp}(\Lambda)=Z_{cell}(\Lambda)-Z_{cell}(0)
\end{equation*}

Using the simple identity
\begin{equation}                                                                              
\label{eq:identity}
C_0 \bigl( (\phi_0 -\phi_\Lambda)\cdot 1\bigr)_\LL=\frac{\Lambda}{D}~
\bigl( \phi_\Lambda \cdot \phi_0\bigr)_\LL
\end{equation}
one writes the spectroscopic impedance as
\begin{equation*}
Z_{sp}(\Lambda)=\frac{\Lambda}{D}~ \frac{\bigl( \phi_\Lambda \cdot \phi_0\bigr)_\LL}
{\bigl(\phi_\Lambda \cdot 1 \bigr)_\LL ~\bigl(\phi_0 \cdot 1 \bigr)_\LL}
\end{equation*}
Applying again the identity (\ref{eq:identity}), one finds a more
convenient form:
\begin{equation*}
Z_{sp}(\Lambda)=\frac{1}{\displaystyle \frac{1}{Z(\Lambda)} - \frac{1}{Z_{cell}(0)}}
\end{equation*}

The new function
\begin{equation*}
Z(\Lambda)=\frac{\Lambda}{D}\bigl(T_\Lambda \phi_0^h \cdot \phi_0^h\bigr)_\LL
\end{equation*}
can be called {\it effective impedance}, where 
\begin{equation*}
\phi_0^h(s)=\frac{\phi_0(s)}{\bigl(\phi_0 \cdot 1 \bigr)_\LL}
\end{equation*}
is the normalized flux density towards the perfectly absorbing working
interface $\pa$. Finally, the spectral decomposition of the spreading
operator $T_\Lambda$ on the basis of eigenfunctions $\V_\alpha$ of the
Dirichlet-to-Neumann operator $\M$ leads to the important relation for
the effective impedance:
\begin{equation}                                                                                    
\label{eq:ZLambda}
Z(\Lambda)=\frac{\Lambda}{D}\sum\limits_\alpha \frac{F_\alpha}{1+\Lambda \mu_\alpha}
\hskip 15mm F_\alpha = \bigl(\phi_0^h \cdot \V_\alpha \bigr)_\LL \bigl(\phi_0^h \cdot 
\V_\alpha^* \bigr)_\LL 
\end{equation}

This relation presents the central result of the continuous approach
developed in \cite{Grebenkov06b}. Let us briefly discuss its physical
meaning.  The spectroscopic impedance $Z_{sp}(\Lambda)$ or,
equivalently, the effective impedance $Z(\Lambda)$, is a physical
quantity that characterizes the transport properties of the whole
working interface. More importantly, this quantity can be measured
directly in experiment (e.g., in electrochemistry). On the other hand,
the local transport properties of the working interface are described
by the single physical parameter $\Lambda$, being related to the
membrane permeability $W$, the electrode resistance $r$ or the
catalyst reactivity $K$ (see Section~\ref{sec:PTL}).  Varying the
parameter $\Lambda$, one changes the local transport properties at
each boundary point and, consequently, the whole linear response of
the working interface. At first sight, one may think that an increase
of the local boundary resistance would imply a proportional increase
of the whole boundary resistance, i.e., $Z(\Lambda)\sim \Lambda$.
This reasoning, being true for a planar surface, becomes invalid in a
general case due to geometrical irregularities and related screening
effects.  In fact, an irregular geometry modifies considerably the
linear response of the working interface
\cite{deLevie65,Nyikos85,Halsey87,Halsey91a,Halsey91b,Halsey92,Sapoval87,Sapoval88,Sapoval89}.
The boundary value problem (\ref{eq:problem_C1}--\ref{eq:problem_C3}),
describing Laplacian transport phenomena on the average, allows
formally to study such geometrical influence.  Practically, however,
this is a very difficult problem. In contrast, the continuous approach
provides an efficient tool to carry out these studies both in
theoretical and numerical ways. In particular, the relation
(\ref{eq:ZLambda}) makes explicit the impedance dependence on the
local transport properties (parameter $\Lambda$) and allows one to
identify contributions due to the {\it physics} and due to the {\it
geometry}, originally involved in the problem in a complex manner. In
other words, whatever the physical problem (diffusion across
semi-permeable membranes, heterogeneous catalysis or electric
transport), the geometry enters {\it only} through the spectral
characteristics of the Dirichlet-to-Neumann operator $\M$: its
eigenvalues $\mu_\alpha$ and the spectral components $F_\alpha$ of the
normalized flux density $\phi_0^h(s)$ on the basis of its
eigenfunctions $\V_\alpha(s)$.

The other important meaning of the relation (\ref{eq:ZLambda}) can be
outlined if one considers the inverse problem \cite{Grebenkov}: what
is the most available information that one can retrieve from a
measurement of the spectroscopic impedance of an unknown working
interface? The mathematical response can be given immediately if one
rewrites (\ref{eq:ZLambda}) as Laplace transform of the new function
$\zeta(\lambda)$:
\begin{equation*}
Z(\Lambda)=\frac{1}{D}\int\limits_0^\infty d\lambda ~e^{-\lambda/\Lambda}~\zeta(\lambda),
\hskip 15mm  \zeta(\lambda)\equiv \sum\limits_\alpha F_\alpha e^{-\lambda \mu_\alpha}
\end{equation*}
Under assumption to be able to measure the impedance with an absolute
precision, one can reconstruct the function $\zeta(\lambda)$ and,
consequently, the set of characteristics $\{\mu_\alpha,~ F_\alpha\}$
which may thus be called {\it harmonic geometrical spectrum} of the
working interface. The hierarchical structure of this spectrum for
self-similar boundaries has been recently investigated
\cite{Grebenkov06d}. Many interesting properties of the function
$\zeta(\lambda)$ remain poorly understood.

\vskip 2mm
In the two following subsections, we are going to discuss some aspects
of the semi-continuous and the discrete descriptions of Laplacian
transport phenomena. These approaches are based on more intuitive
notion of partial reflections on the boundary. Since these
descriptions turn out to be approximations to the continuous approach,
we do not present the circumstantial details.

\subsection{ Semi-continuous Approach }                                                   
\label{sec:approaches_Halsey}

Halsey and Leibig gave the first theoretical description of the
electrolytic double layer response with emphasis on electrochemical
applications \cite{Halsey92}. This description, involving the Green
function of the electrolytic cell, can be reformulated in the
following stochastic language (for details, see \cite{Grebenkov}).
For a given domain $\Omega$ with smooth bounded boundary $\pa$, one
considers the Brownian motion started from a point $x\in\Omega$. When
the diffusing particle hits the boundary at some point $s$, two
complementary events may happen:
\begin{itemize}
\item
with probability $\ve$, the Brownian motion is reflected to the
interior bulk point $s+a n(s)$, slightly above the boundary (here
$n(s)$ is the unit normal vector to the boundary at point $s$, $a$ is
a small positive parameter); the Brownian motion continues from this
point;

\item
or, with probability $1-\ve$, the Brownian motion is terminated at
this point $s$ (absorbed on the boundary).

\end{itemize}
This stochastic process is continued until the absorption on the
boundary and can be called {\it Brownian motion reflected with jump}.
Two new parameters, the jump distance $a$ and the reflection
probability $\ve$, are related to the given physical length
$\Lambda$ \cite{Grebenkov03}:
\begin{equation}                                                                                 
\label{eq:ve}
\ve =\frac{1}{1+ (a/\Lambda)}
\end{equation}

Now, one can calculate the probability $\omega^{(a)}_{x,\Lambda}(s)ds$
that this process is finally absorbed in $ds$ vicinity of the boundary
point $s$. Since the motions before and after each reflection are
independent, this probability can be obtained as the sum of
probabilities to be absorbed after $0$, $1$, $2$, ... reflections:
\begin{equation*}
\begin{split}
\omega^{(a)}_{x,\Lambda}(s)ds & = \bigl[\omega_x(s)ds\bigr] ~(1-\ve ) + \int\limits_\pa 
\bigl[\omega_x(s_1)ds_1\bigr] ~\ve ~ \bigl[\omega_{s_1+an(s_1)}(s)ds\bigr]~ (1-\ve ) + \\
& +\int\limits_\pa  \int\limits_\pa \bigl[\omega_x(s_1)ds_1\bigr] ~\ve ~ 
\bigl[\omega_{s_1+an(s_1)}(s_2)ds_2\bigr]~ \ve ~ \bigl[\omega_{s_2+an(s_2)}(s)ds\bigr]~ 
(1-\ve ) + ~ ... \\
\end{split}
\end{equation*}
For example, the third term represents the probability to hit the
boundary in $ds_1$ vicinity of the point $s_1$, to be reflected to the
neighboring point $s_1+a n(s_1)$, to hit again the boundary in $ds_2$
vicinity of the point $s_2$, to be reflected to the neighboring point
$s_2+a n(s_2)$, to hit the boundary for the last time in $ds$ vicinity
of the point $s$, and to be finally absorbed.  Introducing the
integral operator $Q^{(a)}$, acting from $L^2(\pa)$ to $L^2(\pa)$ as
\begin{equation*}
[Q^{(a)} f](s)=\int\limits_\pa ds'~ f(s') ~\omega_{s'+a n(s')}(s)
\end{equation*}
one rewrites the previous sum as the application of the new integral
operator $T_\Lambda^{(a)}$ to the harmonic measure density
$\omega_x(s)$:
\begin{equation}                                                                           
\label{eq:omega_a}
\omega^{(a)}_{x,\Lambda}(s) = [T_\Lambda^{(a)} \omega_x](s)
\hskip 10mm \textrm{with} \hskip 5mm 
T_\Lambda^{(a)}=(1-\ve)\sum\limits_{k=0}^\infty \bigl(\ve Q^{(a)}\bigr)^k
\end{equation}

What happens when the jump distance $a$ goes to $0$? Hitting the
boundary, the Brownian motion will be reflected to interior points
lying closer and closer to the boundary, i.e., displacements of the
Brownian motion between two serial hits are getting smaller and
smaller. At the same time, the reflection probability $\ve$ tends to
$1$ according to relation (\ref{eq:ve}), i.e., the average number of
reflections increases. Indeed, the distribution of the random number
$\N$ of reflections until the final absorption is simply
\begin{equation}                                                                           
\label{eq:N_geom}
\Pr\{ \N=n \}=(1-\ve)~\ve ^n
\end{equation}
implying that the average number $\E\{ \N \} = \ve (1-\ve)^{-1}$ goes
to infinity.  Does a limiting process exist? The situation is
complicated by the {\it local} choice between reflection and
absorption: at each hitting point, the motion can be absorbed with
vanishing probability $1-\ve$. In order to overcome this difficulty,
one can consider this process from a slightly different point of
view. Actually, one can replace the local condition of the absorption
(with probability $1-\ve$) by its {\it global} analog: the process is
absorbed on the boundary when the number of reflections exceeds a
random variable $\N$ distributed according to the geometrical law
(\ref{eq:N_geom}). Evidently, this modification does not change at all
the properties of the process. At the same time, we gain that the
condition of the absorption becomes independent of the Brownian motion
between serial hits. As a consequence, one can consider the
corresponding limits (as $a\to 0$) separately. So, the Brownian motion
reflected with jump should tend to the reflected Brownian motion as
the jump distance $a$ vanishes. This motion, however, is conditioned
to stop when the number of reflections exceeds the random variable
$\N$. Since the average number $\E\{ \N \}$ goes to infinity in the
limit $a\to 0$, it is convenient to consider a normalized variable
$\chi =a\N$ obeying the following distribution:
\begin{equation}                                                                                     
\label{eq:expon2}     
\Pr\{ \chi \geq \lambda \} = \Pr\{ \N \geq \lambda/a \} = \sum\limits_{[\lambda/a]}^\infty
\Pr\{ \N=n \} \simeq \ve ^{[\lambda/a]} \simeq \exp[-\lambda/\Lambda]
\end{equation}
(the last equality is written with the help of (\ref{eq:ve}) for $a$
going to $0$).  Since the number of reflections on jump distance $a$,
multiplied by $a$, tends to the local time process according to
L\'evy's formula (\ref{eq:Levy2}), the previous condition of
absorption can be reformulated: the motion is absorbed when its local
time process exceeds a random variable distributed according to the
exponential law (\ref{eq:expon2}). One thus concludes that the
Brownian motion reflected with jump should tend to the partially
reflected Brownian motion defined in~\ref{th:PRBM}.

The above analysis, presented as a sketch (without proofs), does not
pretend to a mathematical rigour. It may be considered rather as a
possible justification which can be brought for the semi-continuous
approach if necessary. In particular, one can demonstrate that the
density $\omega^{(a)}_{x,\Lambda}(s)$, given by relation
(\ref{eq:omega_a}), tends to the spread harmonic measure density as
the jump distance $a$ goes to $0$:
\begin{equation*}
\omega_{x,\Lambda}(s)=\lim\limits_{a\to 0} \omega^{(a)}_{x,\Lambda}(s)
\end{equation*}
This relation may be useful for numerical computations (in particular,
it was applied in \cite{Grebenkov06c}). Similarly, the integral
operator $T_\Lambda^{(a)}$ should converge to the spreading operator
$T_\Lambda$ as $a\to 0$. Calculating the geometrical series in
(\ref{eq:omega_a}) and representing $T_\Lambda^{(a)}$ as
\begin{equation*}
T_\Lambda^{(a)}=(1-\ve)\bigl(I-\ve Q^{(a)}\bigr)^{-1}=\left(I+\Lambda ~\frac{I-Q^{(a)}}{a}\right)^{-1}
\end{equation*}
one obtains the following approximation for the Dirichlet-to-Neumann
operator:
\begin{equation*}
\M = \lim\limits_{a\to 0} \frac{I-Q^{(a)}}{a}
\end{equation*}
Again, this relation may be useful for the numerical computation of
this operator.

\vskip 2mm
The advantages of the semi-continuous approach are based on an
apparent intuitive meaning of partial reflections on the
boundary. Moreover, this approach provides even a more realistic
description of physico-chemical processes at microscopic level.  For
example, if one considers the partially reflected Brownian motion
started from a boundary point, the number of hits of the boundary is
infinite for any moment $t>0$ that sounds impossible for real physical
species. The keypoint is that, for diffusion across a semi-permeable
membrane or heterogeneous reaction on a catalytic surface, the
description by the boundary value problem
(\ref{eq:problem_C1}--\ref{eq:problem_C3}) cannot be justified on
length scales less than the mean free path of diffusing particles.
Since the continuous limit $a\to 0$ requires such {\it non-physical}
scales, it is not surprising that the limiting process (the PRBM)
presents some irrealistic properties from the physical point of
view. The similar limitation happens for the electric transport
problem for which the smallest physical scale is given by the
thickness of the double layer, being close to the Debye-H\"uckel
length \cite{Halsey87,Halsey92}.  Evidently, this remark does not
devaluate the efficiency of the continuous approach based on the
partially reflected Brownian motion. On the contrary, the mathematical
rigour of this approach justifies the semi-continuous description and
simplifies its study by introducing the Dirichlet-to-Neumann operator.
However, when dealing with a mathematical description of a physical
problem, one should take care that deduced consequences do not go
beyond the ranges of the model.

The capabilities of the semi-continuous approach are essentially
limited by the fact that the governing operator $Q^{(a)}$ is not
self-adjoint (the function $\omega_{s+an(s)}(s')$ is not symmetric
with respect to the permutation of $s$ and $s'$ except specific
cases). As a consequence, one cannot develop the spectral
decomposition (\ref{eq:ZLambda}) of the impedance. In particular,
there is no possibility to distinguish contributions from different
eigenmodes. Although the operator $Q^{(a)}$ is defined naturally by
the harmonic measure density, it does not provide a proper description
of the problem as it was done with the Dirichlet-to-Neumann operator.

\subsection{ Discrete Approach }                                                             
\label{sec:approaches_Filoche}

Another stochastic approach to Laplacian transport phenomena was
developed by Filoche and Sapoval \cite{Filoche99}. The main idea is to
model the partially reflected Brownian motion by lattice random walks
with partial reflections on the boundary. Actually, one discretizes a
given domain $\Omega$ by $d$-dimensional hypercubic lattice of mesh
$a$ and considers the following stochastic process: started from a
remote source, a random walker jumps to a neighboring site at each
step with probability $(2d)^{-1}$.  When the walker arrives to a
boundary site, it can be reflected to its neighboring site (belonging
to the bulk) with probability $\ve$ (and the motion continues), or it
can be absorbed with probability $(1-\ve)$.  The motion continues
until the final absorption on the boundary, or the return to the
source. One can show \cite{Grebenkov03} that the discrete parameters
$a$ and $\ve$ are related by the expression (\ref{eq:ve}) involving
the continuous physical parameter $\Lambda$.

In the discrete description, the harmonic measure density is replaced
by the distribution of hitting probabilities $(\P_0)_j$ on boundary
sites $j$, (simple) random walks being started from a remote
source. Let $Q^{(a)}_{j,k}$ denote the probability to arrive to the
boundary site $k$ starting from the boundary sites by a random walk in
the bulk without hitting the boundary or the source during the walk%
\footnote{
We use the same notation $Q^{(a)}$ for the integral operator in
semi-continuous approach and for the matrix of these probabilities
since they have the same meaning and even may be used to approximate
each other.}.
%footnote
One can thus calculate the distribution of probabilities
$(\P_\Lambda)_j$ to be finally absorbed on the boundary sites $j$,
when random walks with partial reflections are started from a remote
source. Indeed, the Markov property of this process allows to
calculate $(\P_\Lambda)_j$ as the sum of contributions provided by
random trajectories with $0$, $1$, $2$, ... reflections before the
final absorption:
\begin{equation*}
(\P_\Lambda)_j = (\P_0)_j (1-\ve) + \sum\limits_{k_1} (\P_0)_{k_1} \ve  Q^{(a)}_{k_1,j} 
(1-\ve) + \sum\limits_{k_1}\sum\limits_{k_2} (\P_0)_{k_1} \ve  Q^{(a)}_{k_1,k_2} \ve  
Q^{(a)}_{k_2,j} (1-\ve) + ...
\end{equation*}
(we remind that $\ve=(1+a/\Lambda)^{-1}$).  For example, the second
term represents the product of the following probabilities: to hit a
boundary site $k_1$, to be reflected to its neighboring site, to
arrive to the boundary site $j$, and to be finally absorbed on it. If
one considers $\P_0$ and $\P_\Lambda$ as vectors and $Q^{(a)}$ as
matrix, the summation over intermediate sites $k_1$, $k_2$, ... can be
understood as matrix product:
\begin{equation*}
\P_\Lambda =\left[(1-\ve )\sum\limits_{n=0}^\infty \bigl(\ve Q^{(a)}\bigr)^n \right] \P_0
\end{equation*}
i.e., the distribution of absorption probabilities $(\P_\Lambda)_j$ is
obtained as the application of a linear operator, depending on
$Q^{(a)}$ and $\Lambda$ (or $\ve$), to the distribution of hitting
probabilities $(\P_0)_j$. The symmetric matrix $Q^{(a)}$ represents a
self-adjoint operator, called {\it Brownian self-transport operator}.
Using the normalization property $|Q^{(a)}|\leq 1$ and relation
(\ref{eq:ve}) between $\Lambda $ and $\ve$, one obtains:
\begin{equation*}
\P_\Lambda = T_\Lambda^{(a)} \P_0  \hskip 15mm  T_\Lambda^{(a)}=\left[I+\Lambda~ 
\frac{I-Q^{(a)}}{a}\right]^{-1}
\end{equation*}
The operator $T_\Lambda^{(a)}$, depending on the lattice parameter
$a$, is called {\it (discrete) spreading operator}. The previous
relation, written explicitly as
\begin{equation*}
(\P_\Lambda)_j=\sum\limits_k (\P_0)_k \bigl(T_\Lambda^{(a)}\bigr)_{k,j}
\end{equation*}
allows one to separate random trajectories in two independent parts: 
\begin{itemize}
\item
the random walker started from a remote source arrives to the boundary
site $k$ (first factor);

\item
it continues the motion with partial reflections until the final
absorption on the boundary site $j$ (second factor).
\end{itemize}

One concludes that the absorption probabilities $(\P_\Lambda)_j$
provide a discrete analog of the spread harmonic measure density,
while the matrix $\bigl(T_\Lambda^{(a)}\bigr)_{k,j}$ is a discrete
analog of the kernel $T_\Lambda(s,s')$ of the spreading operator
$T_\Lambda=[I+\Lambda \M]^{-1}$.  In particular, the bounded operators
$(I-Q^{(a)})/a$ can be understood as discrete approximations of the
Dirichlet-to-Neumann operator $\M$ (in resolvent sense).  As for the
semi-continuous approach, we do not furnish the corresponding proofs
(see \cite{Grebenkov} for more details).

\vskip 2mm
The advantage of the discrete description with respect to the
semi-continuous approach is based on the fact that the Brownian
self-transport operator $Q^{(a)}$ and, consequently, the (discrete)
spreading operator $T_\Lambda^{(a)}$ are self-adjoint. This property
allows to employ all the machinery of the spectral theory in order to
express the physical characteristics of Laplacian transport through
eigenmodes of this operator in an explicit way. For example, the
spectral decomposition (\ref{eq:ZLambda}) can be written in the
discrete case. Such decompositions have been used to study Laplacian
transport towards irregular geometries \cite{Grebenkov03,Grebenkov}.
Moreover, the discrete description suggests at least two different
ways to study the problem numerically: direct Monte Carlo simulations
and discrete boundary elements method.

The discrete description, being intuitively the most simple and
useful, may lead to mathematical difficulties when one tries to
proceed the continuous limit $a$ going to $0$. Although the partially
reflected Brownian motion is the natural limit of random walks with
partial reflections, its rigorous demonstration, in our knowledge, is
not yet realized in details.  The interested reader can find more
information on this topic in
\cite{Sabelfeld,Sabelfeld2,Milshtein,Milshtein95,Strook71,Costantini98}.

\section{ Conclusion }

The application of stochastic processes to represent the solution of
boundary value problems is well known and wide used. In particular,
Monte Carlo simulations are generally based on this concept. In this
paper, we gave a brief overview of Laplacian transport phenomena in
different scientific domains (e.g., physics, electrochemistry,
chemistry, physiology) and related stochastic approaches to describe
them. The most attention has been paid to the recently developed
continuous approach based on the partially reflected Brownian motion.
This stochastic process can be thought as rigorous mathematical
description for random trajectories of diffusing particles hitting a
semi-permeable interface, in comparison with more intuitive physical
descriptions by semi-continuous and discrete approaches. The partially
reflected Brownian motion turns out to be the natural limit of the
Brownian motion reflected with jump (semi-continuous approach) and of
the lattice random walks with partial reflections (discrete approach).

The profound relation between the partially reflected Brownian motion
and the spectral properties of the Dirichlet-to-Neumann operator $\M$
are shown to be useful for practical purposes. In particular, the
kernel of the resolvent operator $T_\Lambda=[I+\Lambda \M]^{-1}$ gives
the probability density $T_\Lambda(s,s')$ allowing to reconstruct the
spread harmonic measure $\omega_{x,\Lambda}$.  Moreover, the spectral
decomposition on the complete basis of the Dirichlet-to-Neumann
operator eigenfunctions leads to the explicit analytical formula for
its density. Consequently, the use of the operator $\M$ is an
efficient way to study different probability distributions related to
the partially reflected Brownian motion.

The spectral decomposition of the spectroscopic impedance,
characterizing the linear response of the whole working interface,
leads to an explicit analytical dependence on the physical parameter
$\Lambda$ allowing to identify {\it physical} and {\it geometrical}
contributions which were involved in a complex manner. The harmonic
geometrical spectrum of the working interface contains the complete
information about its transport properties.  The combined use of
stochastic characteristics of the partially reflected Brownian motion
and spectral properties of the Dirichlet-to-Neumann operator opens
encouraging possibilities for further understanding various physical
and chemical transport processes in nature. In this light, a more
profound mathematical analysis of these objects seems to be an
important perspective for the present study.

\subsection*{Acknowledgement}

The author thanks Professor B.~Sapoval and Professor M.~Filoche for
valuable discussions and fruitful collective work on physical aspects
of Laplacian transport phenomena.

\end{document}